\documentclass[12pt,a4paper,twoside]{article}

\newcommand{\pageformat}[6]{\setlength{\hoffset}{-1in}
                  \setlength{\voffset}{-1in}
                  \addtolength{\hoffset}{#5}
                            \addtolength{\voffset}{#6}
                            \setlength{\oddsidemargin}{#1}
                            \setlength{\evensidemargin}{#2}
                            \setlength{\textwidth}{\paperwidth}
                  \addtolength{\textwidth}{-\oddsidemargin}
                  \addtolength{\textwidth}{-\evensidemargin}
                  \addtolength{\textwidth}{-\marginparsep}
                  \addtolength{\textwidth}{-\marginparwidth}
                            \setlength{\topmargin}{#3}
                            \setlength{\textheight}{\paperheight}
                  \addtolength{\textheight}{-\topmargin}
                  \addtolength{\textheight}{-\headheight}
                  \addtolength{\textheight}{-\headsep}
                  \addtolength{\textheight}{-\footskip}
                  \addtolength{\textheight}{-#4}}
\pageformat{2cm}{3cm}{25mm}{25mm}{1pt}{0pt}

\usepackage{ifthen}
\newboolean{article}
    \setboolean{article}{true}
\newboolean{report}
\newboolean{book}
\newboolean{letter}
\newboolean{german}
\newboolean{italian}
\newboolean{nobaselinestretch}
\newboolean{nosectionappendix}
\newboolean{oldtoc}
\newboolean{nosectionequn}
\newboolean{notheorem}

\ifthenelse{\boolean{german}}{
    \usepackage{german}}{}

\usepackage[latin1]{inputenc}

\usepackage{amsmath}
\usepackage{amssymb}
\usepackage[mathscr]{eucal}

\ifthenelse{\boolean{notheorem}}{}{
    \usepackage{theorem}}

\ifthenelse{\boolean{nobaselinestretch}}{}{
    \renewcommand{\baselinestretch}{1.25}}

\newenvironment{env}[2]{\begin{#1}#2\end{#1}}{}
    \newcommand{\beq}[1]{\begin{env}{equation}{#1}}
    \newcommand{\beqn}[1]{\begin{env}{equation*}{#1}}
    \newcommand{\bal}[1]{\begin{env}{align}{#1}}
    \newcommand{\baln}[1]{\begin{env}{align*}{#1}}
    \newcommand{\bga}[1]{\begin{env}{gather}{#1}}
    \newcommand{\bgan}[1]{\begin{env}{gather*}{#1}}
    \newcommand{\bflal}[1]{\begin{env}{flalign}{#1}}
    \newcommand{\bflaln}[1]{\begin{env}{flalign*}{#1}}
    \newcommand{\bmu}[1]{\begin{env}{multline}{#1}}
    \newcommand{\bmun}[1]{\begin{env}{multline*}{#1}}
    \newcommand{\bsp}[1]{\begin{env}{split}{#1}}

    \newcommand{\eeq}{\end{env}}
    \newcommand{\eeqn}{\end{env}}
    \newcommand{\eal}{\end{env}}
    \newcommand{\ealn}{\end{env}}
    \newcommand{\ega}{\end{env}}
    \newcommand{\egan}{\end{env}}
    \newcommand{\eflal}{\end{env}}
    \newcommand{\eflaln}{\end{env}}
    \newcommand{\emu}{\end{env}}
    \newcommand{\emun}{\end{env}}
    \newcommand{\esp}{\end{env}}

\newcommand{\lf}{\vspace{2ex}}

\renewcommand{\bf}[1]{\textbf{#1}}
\renewcommand{\it}[1]{\textit{#1}}

\renewcommand{\sf}[1]{\textsf{#1}}

\renewcommand{\tt}[1]{\texttt{#1}}
\newcommand{\hl}[1]{\bf{\it{#1}}}
\newcommand{\mrm}[1]{\mathrm{#1}}
\newcommand{\mbf}[1]{\mathbf{#1}}
\newcommand{\msf}[1]{\text{\small $\sf{#1}$}}

\newcommand{\cmc}[1]{\mathcal{#1}}
\newcommand{\eus}[1]{\mathscr{#1}}

\newcommand{\bb}[1]{\mathbb{#1}}

\newcommand{\mscriptsize}[1]{{\setlength{\arraycolsep}{.3ex}\text{\scriptsize$#1$}}}
\newcommand{\mtiny}[1]{{\setlength{\arraycolsep}{.3ex}\text{\tiny$#1$}}}
\newcommand{\nbd}[1]{$#1$\nobreakdash--}
\newcommand{\ol}[1]{\overline{#1}}

\newcommand{\vt}{\vartheta}

\newcommand{\vp}{\varphi}

\newcommand{\abs}[1]{\left\lvert#1\right\rvert}
\newcommand{\norm}[1]{\left\lVert#1\right\rVert}

\newcommand{\family}[1]{\left(#1\right)}

\newcommand{\bfam}[1]{\bigl(#1\bigr)}

\newcommand{\AB}[1]{\langle#1\rangle}

\newcommand{\CB}[1]{\{#1\}}
\newcommand{\bCB}[1]{\bigl\{#1\bigr\}}
\newcommand{\BCB}[1]{\Bigl\{#1\Bigr\}}

\newcommand{\Matrix}[1]{\begin{pmatrix}#1\end{pmatrix}}

\newcommand{\sMatrix}[1]{\mscriptsize{\Matrix{#1}}}
\newcommand{\rsMatrix}[1]{\raisebox{.18ex}{\sMatrix{#1}}}
\newcommand{\tMatrix}[1]{\mtiny{\Matrix{#1}}}
\newcommand{\rtMatrix}[1]{\raisebox{.3ex}{\tMatrix{#1}}}

\newcommand{\set}[2][]{
    \ifthenelse{\equal{#1}{}}{
        \CB{#2}}{
        \CB{#1~|~#2}}}
\newcommand{\bset}[2][]{
    \ifthenelse{\equal{#1}{}}{
        \bCB{#2}}{
        \bCB{#1~|~#2}}}
\newcommand{\Bset}[2][]{
    \ifthenelse{\equal{#1}{}}{
        \BCB{#2}}{
        \BCB{#1~\big|~#2}}}
\newcommand{\zero}{\CB{0}}

\DeclareMathOperator{\ls}{\normalfont\msf{span}}

\DeclareMathOperator{\ssls}{\scriptstyle\sf{span}}
\DeclareMathOperator{\Sls}{\text{\scriptsize $\sf{span}$\,}}
\DeclareMathOperator{\cls}{\ol{\ls}}

\DeclareMathOperator{\sscls}{\ol{\ssls}}
\DeclareMathOperator{\Scls}{\ol{\Sls}}

\DeclareMathOperator{\id}{\normalfont\msf{id}}

\renewcommand{\ker}{\operatorname{\msf{ker}}}

\newcommand{\C}{\bb{C}}

\newcommand{\E}{\bb{E}}

\newcommand{\K}{\bb{K}}

\newcommand{\Q}{\bb{Q}}

\newcommand{\T}{\bb{T}}

\newcommand{\cA}{\cmc{A}}
\newcommand{\cB}{\cmc{B}}
\newcommand{\cC}{\cmc{C}}
\newcommand{\cD}{\cmc{D}}

\newcommand{\cI}{\cmc{I}}

\newcommand{\sB}{\eus{B}}

\newcommand{\sK}{\eus{K}}

\newcommand{\sT}{\eus{T}}

\newcommand{\U}{\mbf{1}}
\newcommand{\F}{{\mrm{F}}}
\newcommand{\G}{\Gamma}

\ifthenelse{\boolean{nosectionequn}}{}{
    \numberwithin{equation}{section}
    }

\ifthenelse{\boolean{article}\or\boolean{letter}\or\boolean{nosectionequn}}{
    \setboolean{nosectionappendix}{true}}{}
\ifthenelse{\boolean{nosectionappendix}}{}{
    \renewcommand{\appendix}{
        \chapter*{\appendixname}
        \addcontentsline{toc}{chapter}{\appendixname}
        \renewcommand{\thesection}{\Alph{section}}
        \setcounter{section}{0}}}
   
\ifthenelse{\boolean{report}\or\boolean{book}}{
    }{}

\ifthenelse{\boolean{notheorem}}{}{
        \newcommand{\notename}{Note.}
        \newcommand{\mnname}{Mathematical note.}
        \newcommand{\enname}{End of the note.}
        \newcommand{\definame}{Definition.}
        \newcommand{\propname}{Proposition.}
        \newcommand{\lemname}{Lemma.}
        \newcommand{\exname}{Example.}
        \newcommand{\exername}{Exercise.}
        \newcommand{\remname}{Remark.}
        \newcommand{\obname}{Observation.}
        \newcommand{\thmname}{Theorem.}
        \newcommand{\corname}{Corollary.}
        \newcommand{\proofname}{Proof.}
    \ifthenelse{\boolean{german}}{
        \renewcommand{\mnname}{Mathematische Notiz.}
        \renewcommand{\enname}{Ende der Notiz.}
        \renewcommand{\exname}{Beispiel.}
        \renewcommand{\exername}{Übung.}
        \renewcommand{\remname}{Bemerkung.}
        \renewcommand{\obname}{Beobachtung.}
        \renewcommand{\thmname}{Satz.}
        \renewcommand{\corname}{Korollar.}
        \renewcommand{\proofname}{Beweis.}}{}
    \ifthenelse{\boolean{italian}}{
        \renewcommand{\mnname}{Nota matematica.}
        \renewcommand{\enname}{Fina della nota.}
        \renewcommand{\definame}{Definizione.}
        \renewcommand{\propname}{Proposizione.}
        \renewcommand{\exname}{Esempio.}
        \renewcommand{\exername}{Esercizio.}
        \renewcommand{\remname}{Nota.}
        \renewcommand{\obname}{Osservazione.}
        \renewcommand{\thmname}{Teorema.}
        \renewcommand{\corname}{Corollario.}
        \renewcommand{\proofname}{Dimostrazione.}

       \renewcommand{\appendixname}{Appendice}

       }{}
    \theoremheaderfont{\normalfont\bfseries}
    \theoremstyle{change}
        \theorembodyfont{\rmfamily}
            \newtheorem{emp}{}[section]
                \newcommand{\bemp}[1][]{
                    \begin{emp}\hskip-\labelsep\bf{#1}\hskip\labelsep}
                \newcommand{\eemp}{\end{emp}}
\newtheorem{itemp}[emp]{}
                \newcommand{\bitemp}[1][]{
                    \begin{itemp}\hskip-\labelsep\bf{#1}\hskip\labelsep\normalfont\itshape}
                \newcommand{\eitemp}{\end{itemp}}
            \newtheorem{note}[emp]{\notename}
                \newcommand{\bnote}{\begin{note}}
                \newcommand{\enote}{\end{note}}
            \newtheorem{mn}[emp]{\mnname}
                \newcommand{\bnm}{\begin{mn}~\begin{quotation}\renewcommand{\baselinestretch}{1}\small\noindent\ignorespaces}
                \newcommand{\enm}{\end{quotation}\hfill\bf{\enname}\end{mn}}
            \newtheorem{ex}[emp]{\exname}
                \newcommand{\bex}{\begin{ex}}
                \newcommand{\eex}{\end{ex}}
            \newtheorem{exer}[emp]{\exername}
                \newcommand{\bexer}{\begin{exer}}
                \newcommand{\eexer}{\end{exer}}
            \newtheorem{defi}[emp]{\definame}
                \newcommand{\bdefi}{\begin{defi}}
                \newcommand{\edefi}{\end{defi}}
            \newtheorem{rem}[emp]{\remname}
                \newcommand{\brem}{\begin{rem}}
                \newcommand{\erem}{\end{rem}}
            \newtheorem{ob}[emp]{\obname}
                \newcommand{\bob}{\begin{ob}}
                \newcommand{\eob}{\end{ob}}
        \theorembodyfont{\normalfont\itshape}
            \newtheorem{thm}[emp]{\thmname}
                \newcommand{\bthm}{\begin{thm}}
                \newcommand{\ethm}{\end{thm}}
            \newtheorem{prop}[emp]{\propname}
                \newcommand{\bprop}{\begin{prop}}
                \newcommand{\eprop}{\end{prop}}
            \newtheorem{cor}[emp]{\corname}
                \newcommand{\bcor}{\begin{cor}}
                \newcommand{\ecor}{\end{cor}}
            \newtheorem{lem}[emp]{\lemname}
                \newcommand{\blem}{\begin{lem}}
                \newcommand{\elem}{\end{lem}}
\newenvironment{empn}[1]{\lf\noindent\bf{#1}\ignorespaces\hskip\labelsep}{\lf}
		\newcommand{\bempn}[1]{\begin{empn}{#1}}
		\newcommand{\eempn}{\end{empn}}
		\newcommand{\bitempn}[1]{\begin{empn}{#1}\normalfont\itshape}
		\newcommand{\eitempn}{\end{empn}}
                \newcommand{\bnmn}{\begin{empn}{\mnname}~\begin{quotation}\renewcommand{\baselinestretch}{1}\small\noindent\ignorespaces}
                \newcommand{\enmn}{\end{quotation}\hfill\bf{\enname}\end{empn}}
		\newcommand{\bexn}{\begin{empn}{\exname}}
		\newcommand{\eexn}{\end{empn}}
		\newcommand{\bexern}{\begin{empn}{\exername}}
		\newcommand{\eexern}{\end{empn}}
		\newcommand{\bdefin}{\begin{empn}{\definame}}
		\newcommand{\edefin}{\end{empn}}
		\newcommand{\bremn}{\begin{empn}{\remname}}
		\newcommand{\eremn}{\end{empn}}
		\newcommand{\bobn}{\begin{empn}{\obname}}
		\newcommand{\eobn}{\end{empn}}

		\newcommand{\bpropn}{\bitempn{\propname}}
		\newcommand{\epropn}{\eitempn}

\newcommand{\qedsymbol}{~\rule[-0.35mm]{2mm}{2mm}}
    \newcounter{proof}[emp]
    \newenvironment{Proof}[1]{
        \vspace{1ex}
        \renewcommand{\item}[1][\stepcounter{proof}(\roman{proof})]%
            {##1\hskip\labelsep}
        \noindent\textsc{#1\hskip\labelsep}}{
        \nolinebreak\qedsymbol}
    \newcommand{\proof}[1][\proofname]{
        \begin{Proof}{#1}\ignorespaces}
    \newcommand{\qed}{\end{Proof}}
    \newcommand{\noqed}{
        \renewcommand{\qedsymbol}{}
        \end{Proof}}}
    \ifthenelse{\boolean{italian}}{
        \renewcommand{\proofname}{Dimostrazione.}}{}

\usepackage[varg]{txfonts}

\usepackage{stmaryrd}

\renewcommand{\thefootnote}{[\arabic{footnote}]}

\usepackage{hyperref}

\newcommand{\sds}{\hspace{.1ex}.\hspace{.2ex}}

\begin{document}

\bibliographystyle{amsalpha}

\title{Ideal Submodules\\\it{versus} Ternary Ideals\\\it{versus} Linking Ideals}

\author{Michael Skeide{\renewcommand{\thefootnote}{}}
\setcounter{footnote}{0}
}

\date{This revision December 2020}

\maketitle 


\begin{abstract}
\noindent
We show that ideal submodules and closed ternary ideals in Hilbert modules are the same. We use this insight as a little peg on which to hang a little note about interrelations with other notions regarding Hilbert modules. In Section \ref{idealSEC}, we show that the ternary ideals (and equivalent notions) merit fully, in terms of homomorphisms and quotients, to be called ideals of (not necessarily full) Hilbert modules. The properties to be checked are intrinsically formulated for the modules (without any reference to the algebra over which they are modules) in terms of their ternary structure. The proofs, instead, are motivated from a third equivalent notion, linking ideals (Section \ref{liSEC}), and a Theorem (Section \ref{idealSEC}) that all extends nicely to (reduced) linking algebras. As an application, in Section \ref{extSEC}, we introduce ternary extensions of Hilbert modules and prove most of the basic properties (some new even for the known notion of extensions of Hilbert modules), by reducing their proof to the well-known analogue  theorems about extensions of \nbd{C^*}algebras. Finally, in Section \ref{qSEC}, we propose several open problems that our method naturally suggests.
\end{abstract}

\setcounter{section}{-1}
\section{Ideal submodules \it{versus} ternary ideals} \label{zero}

Let us start right with our \it{peg} from the abstract: Ideal submodule=closed ternary ideal.

\bempn{Definition ISM.~}
Let $I$ be a closed ideal in a \nbd{C^*}algebra $\cB$, and let $E$ be a Hilbert \nbd{\cB}module. The \hl{ideal submodule} of $E$ associated with $I$ is $\cls EI$. More generally, we say $K\subset E$ is an \hl{ideal submodule} of $E$ if $K=\cls EI$ for some ideal $I$ in $\cB$.
\eempn

\bempn{Definition TI.~}
A linear subspace $K$ of a Hilbert \nbd{\cB}module $E$ is a \hl{ternary ideal} if $E\AB{K,E}\subset K$.
\eempn

\bempn{Definition RI.~}
For a Hilbert \nbd{\cB}module $E$ we denote by $\cB_E:=\cls\AB{E,E}$ its \hl{range ideal}.
\eempn

We use the same notation $\cB_K:=\cls\AB{K,K}$ for any subset $K$ of $E$. But, only if $K$ is a \nbd{\cB}submodule, is it granted that $\cB_K$ is an ideal in $\cB$.

\bempn{Proposition E.~}
For a subset $K$ of a Hilbert \nbd{\cB}module $E$ the following conditions are equivalent:
\begin{enumerate}
\item
$K$ is an ideal submodule of $E$.

\item
$K$ is a closed ternary ideal in $E$.
\end{enumerate}
\eempn

\vspace{-3ex}
\proof
Of course, if $K$ is an ideal submodule, it is a closed \nbd{\cB}submodule of $E$. Like for any closed submodule of $E$, this implies $\cls K\cB_K=K$ and $\cls\AB{K,E}=\cB_K=\cls\AB{K,E}$. From the first property we infer that $\cB_K$ is the unique smallest ideal in $\cB$ for which $K$ is the associated ideal submodule. (Indeed, $K=\cls EI=\cls EI\cB_K$. So, $\cls I\cB_K=I\cap\sB_K$ is a smaller ideal with which $K$ is associated. Since $\cB_K=\cls\AB{K,K}=\cls \AB{KI,KI}$, the ideal $I\cap\cB_K$ cannot be smaller than $\cB_K$.) From the second property we infer $\AB{K,E}\subset\cB_K$, so, $E\AB{K,E}\subset E\cB_K\subset K$.

\it{Vice versa} if $K$ is a closed ternary ideal, then
\beqn{
\AB{E,E}\AB{K,E}
~\subset~
\AB{E,K}.
}\eeqn
Making use of an approximate unit for $\cB_E$, we get $\AB{K,E}\subset\cls\AB{E,K}$, and by taking adjoints $\AB{E,K}\subset\cls\AB{K,E}$, so $\cls\AB{K,E}=\cls\AB{E,K}$. Since $\AB{K,E}\AB{E,E}\subset\AB{K,E}$, we see that $I:=\cls\AB{K,E}$ is an ideal in $\cB_E$ and, further, in $\cB$. Clearly, $EI\subset K$. On the other hand, $\cls EI=\cls E\AB{E,K}=\cls\sK(E)K\supset K$. (Use an approximate unit for $\sK(E)$.) In conclusion, $K=\cls EI$ (and, of course, $I=\cls\AB{K,E}=\cB_K$).\qed%
\footnote{
This contradicts \cite[Theorem 4.3]{Kol17}. Indeed, the example proposed in its proof has the general structure $E=\cB\oplus\cB$ and $K=I_1\oplus I_2$ for distinct closed ideals $I_1\ne I_2$ of $\cB$. (Direct sum of Hilbert \nbd{\cB}modules. In the example we even have $I_1\supset I_2$.) Since $\Scls EI=I\oplus I$ for all $I$, there is no $I$ such that $\Scls EI=K$, so $K$ is not an ideal submodule. But, also $\Scls E\AB{K,E}=\Scls E(I_1+I_2)$ (not $\oplus$!) is not contained in $K$. So, unlike what is claimed in \cite{Kol17}, $K$ is not a closed ternary ideal, and the example is not a counter example.
}

\lf\noindent
\bf{Conventions.~}
Here and in the sequel, we assume that notions like \it{Hilbert modules}, $E^*\ni x^*=\AB{x,\bullet}$, \it{rank-one} operators $xy^*\colon z\mapsto x\AB{y,z}$, \it{compact operators} $\sK(E,F)=\cls FE^*$, and their basic properties are known. We do \bf{not} adopt the common standard convention according to which writing a product of spaces would mean the closure of the linear span: $AB$ for subsets of spaces for which a product $ab$ is defined, means exactly the set $AB=\CB{ab\colon a\in A,b\in B}$ and nothing else. $\sB^a(E,F)$ means the set of \it{adjointable} (automatically bounded) operators between Hilbert modules $E$ and $F$. When necessary, we add a superscript $^{bil}$ to mean \it{\nbd{\cA}\nbd{\cB}linear maps} between \nbd{\cA}\nbd{\cB}bimodules.

\newpage

\section{Introduction} \label{intro}

What is our main issue? Well, after having set up the definitions, our \it{peg} (that is, Proposition E in Section \ref{zero}) is, actually, not much more than an elaborate triviality -- a property it shares with many more statements ranging from very well-known to new to be found in the literature and later on in these notes: Once you have the statement you wish to prove, the proofs are easy/straightforward/to be left as simple exercises. So, why contribute to the long list of existing papers (which we largely omit) another one?

Towards answering this question, let us have a closer look at our \it{peg}. There are (at least) two, good (we think), reasons why  the reader and we can benefit from elaboration of our \it{peg}.
\begin{enumerate}
\item
\it{Working:}
If there are two equivalent definitions of the same thing, \it{closed ideals} in Hilbert modules, then it is natural to ask if one of them is easier to work with. This, actually, splits further into two subquestions. Firstly, \it{easier to work with} in the sense of fewer or easier conditions to be checked to apply the results. Secondly, \it{easier to work with} in the sense of finding proofs more easily. A good deal of these notes is dedicated to illustrate that our choices, closed \it{ternary ideals} and other ternary structures (quotients, morphisms, extensions ...), appear \it{unbeatable} when it comes to verify conditions. As for proving results, we soon find it convenient to introduce a third (of course, equivalent) notion, closed \it{linking ideals} (Section \ref{liSEC}), based on the well-known fact that every Hilbert module embeds into its (reduced) \it{linking algebra}.

\item
\it{Motivation:}
None of the two definitions of \it{ideals} in Hilbert modules appears actually very well motivated. While both appear reasonable, also demonstrated by the results that can be proved for them, they still appear rather \it{ad hoc} in the style ``let us try a definition of ideals in Hilbert modules (without telling why) and see what we can prove for them''. In Section \ref{idealSEC}, we show that closed \it{ternary ideals} in Hilbert modules do exactly what we expect by analogy with closed \it{ideals} in \nbd{C^*}algebras, while \it{ideal submodules} do not (unless we limit our attention to full Hilbert modules).

The motivation via linking ideals works even better. The closed ideals of the reduced linking algebra of a Hilbert module are exactly the linking algebras of its closed ternary ideals. Blockwise homomorphisms correspond to ternary homomorphisms, quotients to quotients. New notions are motivated by their known analogues for the linking algebras, results about them follow from well-know results applied to the linking algebras. In Sections \ref{extSEC} we apply this to \it{extensions} of Hilbert modules, while in Section \ref{qSEC} we propose a number of new notations and problems we do not resolve in these notes.
\end{enumerate}

\noindent
As we said above, most results have simple proofs, and many of them (though, surely not all) are known. It is our intention to underline this simplicity. Necessarily, this scope requires the exposition to be (on a consistent level) self-contained. It is surely not a good idea to send the reader to existent proofs. (Partly, not to interrupt the flow; partly, because the arguments we exploit in the proofs here are \it{multi-use}; and partly, because the existing proofs frequently appear as part of a different, usually larger, context, and make reference to this context. (For instance, the book \cite{BlLMe04} by Blecher and Le Merdy is an excellent reference for many of the known notions and their properties; however, their discussion is embedded in the context of \it{operator spaces} and not entirely independent of that context.) Also the simplicity of the (known) results leads to the fact that they occur in many places in the literature, frequently even without knowing of each other's existence. It is, therefore, often difficult, if not impossible, to figure out the primary source. Additionally, some of the results and notions do occur, yes, but in a different \it{clothing} (for instance, the \it{Rieffel correspondence}) which, when explained the moment they occur here, would sensibly disturb the flow and the simplicity of the argument. We do make a certain effort to provide references at least for some of the best known results. But in order to not obstruct the simplicity, we will mostly do this in the form of separated remarks towards the sections' ends or, occasionally, in footnotes.

Let us just give an instance of the difficulties, when wishing to do full justice to earlier work. Apart from our Definition TI of \it{ternary ideal}, there exists already a definition in the context of  \it{ternary operator spaces} (TRO) (see also Remark \ref{TROrem} and Observation \ref{TROob}). After Proposition E in Section \ref{zero}, our second new result is equivalence of the two definitions of \it{ternary ideal} as stated in Corollary \ref{TI=TIcor}. A third related equivalent substructure, \it{subbimodules} of so-called \it{Hilbert bimodules} will be illustrated in the Remark \ref{MERCrem} about \it{Morita equivalence} and \it{Rieffel correspondence}. The scope of that remark is taking the time to sort out a quite enormous \it{salad} of similar but not identical notions. It also illustrates how difficult it actually is not only to find existing results, but also to carefully distinguish between existing results and ``almost existing'' results. (Remark \ref{MERCrem} and other remarks in these notes do not help to understand these notes; they are just an attempt to be as \it{fair and square} as possible without having to write a little comparative PhD-thesis. Mathematically, it is completely safe for any reader to skip them.)

Apart from many little side results that range in between new and ``almost existing'' (and which we do not list here), we like to mention in particular the (apparently) new results Theorem \ref{qthm} and the \it{addition of extensions} in \ref{Badd}. But, also the notion of \it{ternary extension} (Definition \ref{TEdefi}) appears new. Hence, the whole Section \ref{extSEC}, in a sense, is new. After a reinterpreting using the equivalences established in Sections \ref{liSEC} and \ref{idealSEC}, the results (apart from \ref{Badd}) become equivalent to (parts of) the results by Bakic and Guljas \cite{BaGu04,BaGu03}. However, using the equivalences established in Sections \ref{liSEC} and \ref{idealSEC} in a different way, we obtain these result by reduction to well-known results about extensions of \nbd{C^*}algebras applied to the (reduced) linking algebras, without having to redo any proof as in \cite{BaGu04,BaGu03}. In Section \ref{qSEC}, we propose the (apparently) new notions of \it{conditional expectations}, \it{semisplit extensions}, and \it{linking hereditary subspaces} for Hilbert modules, and state open problems for them. We compare \it{linking hereditary subspaces} with \it{ternary hereditary subspaces} (the latter also known as \it{inner ideals} of TROs), and derive some preliminary results about them. The question if these two notions are equivalent, we leave open.

\newpage
For the balance of this introduction, we already address to some extent our answer to Point 1 above. Proposition E tells us that (modulo closure) Definition ISM and Definition TI distinguish the same collection of subspaces of a Hilbert module. So what are the differences? To examine whether a given subset $K$ of a Hilbert \nbd{\cB}module $E$ is an ideal submodule, we have to examine whether there exists (or not) an ideal $I$ in $\cB$ such that the elements in $EI$ are in $K$ and, there, total. (Well, to be honest, the original definition in Bakic and Guljas \cite[Definition 1.1]{BaGu02} speaks only of the ideal submodule associated with a given ideal $I\subset\cB$; we added already the part with ``for some $I$'', directing attention towards the question when a given $K$ is an ideal submodule.%
\footnote{
It is noteworthy, that the ideal submodule of $E$ associated with $I$ can also be described as $\CB{x\in E\colon\AB{x,x}\in I}$. The simple proof is in Takahashi \cite[Lemma 2.07]{TakA79}, and in this form \cite{TakA79} is the first appearance of ideal submodules we could get hold of. (Except for dealing with Hilbert \it{left} modules (see also Footnote \ref{MEsymFN}), the paper is nicely written and a recommended reading. But some of its results were definitely not new at that time. See also Remark \ref{MERCrem}.)
}%
) Apart from the problem of having to guess what (a suitable) $I$ is (or even  to exclude its existence), this is quite a number of conditions to be verified. On the contrary, for checking whether a subspace $K$ is a ternary ideal, there is just the ternary invariance condition $E\AB{K,E}\subset K$ to be controlled. (Well, to be honest, ideal submodules are (linear) subspaces for free while for ternary ideals we have to require that.)

So, the notion of ideal submodule of $E$ puts emphasis on the ideal $I\subset\cB$ and, therefore, on the \nbd{C^*}algebra $\cB$ over which $E$ is a Hilbert module; this goes through the whole paper \cite{BaGu02} (and its successors \cite{BaGu04,BaGu03}) in the sense that the homomorphisms $v$ between modules over \nbd{C^*}algebras are maps for which there have to exist (and to be identified!) homomorphisms $\vp$ between the \nbd{C^*}algebras fulfilling a condition (see \eqref{phiiso}) relating the two maps $v$ and $\vp$. When working with these notions, we first have to identify something (an ideal in the \nbd{C^*}algebra or a homomorphism between the \nbd{C^*}algebras) that relates to the \nbd{C^*}algebra(s); and the definitions do not indicate how to find these somethings. On the contrary, the condition to be a ternary ideal relies only on the ternary product $(x,y,z)\mapsto x\AB{y,z}$; it is, therefore, intrinsic to the \it{ternary operator space} structure that uniquely characterizes a Hilbert module without (explicit) reference to the \nbd{C^*}algebra over which it is a module. The same remains true for the homomorphisms $v$, which will be \it{ternary} (see Section \ref{idealSEC}).

Of course, there is a good candidate for $I$, namely $\cB_K=\cls\AB{K,K}$. (Finding $\vp$ for $v$ is harder and, in fact, not always possible if the modules are not required to be full; see Section \ref{liSEC}.) But nothing in Definition ISM suggests this; the insight that $\cB_K$ is a good candidate is, actually, already part of the proof of Proposition E. Also, if $K$ is an ideal submodule, the ideal $I$ illustrating this is not unique. (We leave it as a little exercise to deduce from the right line in the proof of Proposition E that $\cls EI_1=\cls EI_2$ if and only if $I_1\cap\cB_E=I_2\cap\cB_E$, and invite the reader to find $E$ and $I_1\ne I_2$ satisfying this condition.)

Of course, if $K\AB{E,E}\subset K$, then $\cB_K$ is an ideal in $\cB_E$, hence $\cB$. (Effectively, $\ol{K}$ is a closed submodule of $E$ if and only if $K\AB{E,E}\subset\ol{K}$.) The following little discussion shows (Corollary \ref{TI=TIcor}) that is is enough (and necessary) to add the symmetric condition $E\AB{E,K}\subset K$ to guarantee that $K$ is a ternary ideal:

Looking at the proof of Proposition E, in either direction it was important to note that (for different reasons in each case) $\cls\AB{K,E}=\cls\AB{E,K}$. This is true for all submodules of $E$, so we get this as soon as we established that $K$ is a submodule. While in the forward direction this was obvious, for the backward direction we used some case specific trickery. One can also argue, as we do in the following corollary, appealing to the general theorem that a two-sided closed ideal in a \nbd{C^*}algebra is a \nbd{*}ideal (see, for instance, Murphy \cite[Theorem 3.1.3]{Mur90}).

\bcor \label{TI=TIcor}
Let $K$ be a subset of a Hilbert \nbd{\cB}module $E$.
\begin{enumerate}
\item \label{TI1}
If $K$ is a closed ternary ideal, then
\beq{ \label{TI=TI}
E\AB{E,K}\subset K
\text{~~~~~~and~~~~~~}
K\AB{E,E}\subset K.
}\eeq

\item \label{TI2}
Both conditions in \eqref{TI=TI} together imply that $\ol{K}$ is a ternary ideal.

\item \label{TI3}
On the contrary, none of the two conditions alone implies that $\ol{K}$ is a ternary ideal.
\end{enumerate}
\ecor

\vspace{-1ex}
\proof
\eqref{TI1}~~The second condition in \eqref{TI=TI} is true for every submodule $K$. If $K$ is a closed ternary ideal then the first condition follows from $\cls\AB{E,K}=\cls\AB{K,E}$ as in the proof of Proposition E.

\eqref{TI2}~~Conversely, suppose both conditions hold. Then $\cls\AB{E,K}$ is a closed two-sided ideal of $\cB$, hence, a \nbd{*}ideal. Therefore, $\cls\AB{E,K}=\cls\AB{K,E}$, hence, $E\AB{K,E}\subset\cls E\AB{E,K}\subset\ol{K}$.

\eqref{TI3}~~Finally, taking $E=\cB\ni\U$, the subspaces $K$ fulfilling only the first and only the second condition are the left and the right ideals in $\cB$, respectively, while (we just used that) both conditions together mean $K$ is an ideal. Of course, there are left ideals that are not right ideals, and \it{vice versa}.\qed

\brem \label{TROrem}
Hilbert modules are TROs (see Section \ref{liSEC}). Closed subspaces $K$ satisfying the conditions in \eqref{TI=TI} have been called \it{triple ideals} of the TRO $E$ in Blecher and Le Merdy \cite[ahead of Corollary 8.3.4]{BlLMe04} and \it{ternary ideals} of the TRO $E$ in Blecher and Neal \cite[p.230]{BlNe07}. (We shall use the latter one,  Blecher and Neal \cite{BlNe07}, as a basic reference for structures related to TROs. Not only is it nearer to what we need here and has everything tied together in its Section 1, while in \cite{BlLMe04} the arguments are scattered over several chapters. Also the definition of TRO (certain closed subspaces of \nbd{C^*}algebras) we use in our Section \ref{liSEC}, is specifically the one from \cite{BlNe07}, while \cite[3.1.2]{BlLMe04} defines TRO as certain closed subspaces ``of $\sB(H)$, or of $\sB(K,H)$''.)

Our result in Corollary \ref{TI=TI}, saying that \it{ternary ideals} of the Hilbert module $E$ in the sense of \cite{BlNe07} are exactly the \it{closed ternary ideals} in the sense of our Definition TI, appears to be new.
\erem

\newpage

\section{... \it{versus} linking ideals} \label{liSEC}

Let us complete the picture about ideals in Hilbert modules -- and the title of Section \ref{zero} -- by adding a third point of view. Recall that the \hl{linking algebra} of a Hilbert \nbd{\cB}module $E$ is $\sK\rtMatrix{\cB\\E}=\rtMatrix{\sK(\cB)&\sK(E,\cB)\\\sK(\cB,E)&\sK(E)}=\rtMatrix{\cB&E^*\\E&\sK(E)}$, while the \hl{reduced linking algebra} of $E$ is $\sK\rtMatrix{\cB_E\\E}=\rtMatrix{\cB_E&E^*\\E&\sK(E)}$. The two coincide if and only $E$ is \hl{full}, that is, if $\cB_E=\cB$.

\bdefi \label{lidef}
A subset $K$ of a Hilbert \nbd{\cB}module $E$ is a \hl{linking ideal} if $\rtMatrix{\sscls\AB{K,K}&K^*\\K&\sscls KK^*}$ is an ideal in the reduced linking algebra of $E$.
\edefi

\vspace{-2ex}
\bthm \label{Ethm}
For a subset $K$ of a Hilbert \nbd{\cB}module $E$ the following conditions are equivalent:
\begin{enumerate}
\item \label{I1}
$K$ is an ideal submodule of $E$.

\item \label{I2}
$K$ is a closed ternary ideal in $E$.

\item \label{I3}
$K$ is a closed linking ideal in $E$.
\end{enumerate}
\ethm

\lf\noindent
We prefer to discuss a part of the proof in a more general situation. Recall that a \hl{ternary subspace} of a Hilbert module $E$ is a linear subspace $F\subset E$ such that $F\AB{F,F}\subset F$.

\blem \label{TsSlem}
If $F$ is a closed ternary subspace of a Hilbert \nbd{\cB}module $E$, then $\cB_F$ is a \nbd{C^*}sub\-al\-ge\-bra of $\cB$ and $F$ is a (full) Hilbert \nbd{\cB_F}module. The linking algebra of $F$ is a \hl{blockwise} subalgebra of the linking algebra of $E$ (that is, $\sK\rtMatrix{\cB_F\\F}$ is a \nbd{*}subalgebra of $\sK\rtMatrix{\cB\\E}$ with respect to the canonical identifications $\cB_F\subset\cB$, $F^*\subset E^*$, $F\subset E$, and $\sK(F)\subset\sK(E)$).

Moreover, the \nbd{C^*}subalgebra generated by $F$ is
\beqn{
\rsMatrix{\Scls\AB{F,F}&F^*\\F&\Scls FF^*}
}\eeqn
and coincides with the reduced linking algebra $\sK\rtMatrix{\cB_F\\F}$.

$F$ is a Hilbert submodule of $E$ (if and only) if $\cB_F$ is an ideal in $\cB$.

Of course, every Hilbert submodule of $E$ is a ternary subspace.
\elem

\proof
Straightforward verification.\qed

\bob \label{TROob}
(See also Remark \ref{MERCrem}.) The lemma applies, in particular, to a ternary subspace $F$ of a \nbd{C^*}algebra $\cA$, a so-called \hl{ternary ring of operators} (\hl{TRO})%
\footnote{
Note that we do note require that $F$ be closed as in \cite{BlNe07}. Also, it does not matter whether we consider $\cA$ a Hilbert \nbd{\cA}module (to apply our definition of ternary subspace) or as \nbd{C^*}algebra (to apply the one from  \cite{BlNe07}).

Abstract \it{ternary rings} occurred (at least) as early as Hestenes \cite{Hes62}, which is, clearly, motivated by the TRO generated by a single matrix in Hestenes \cite{Hes61}. Zettl \cite{Zet83} introduces abstract \it{ternary \nbd{C^*}rings} and shows how to modify the ternary product of a ternary \nbd{C^*}ring by means of a (unique, self-inverse) isomorphism of ternary rings in order to make it isomorphic to a (closed) TRO.
}
in the definition of Blecher and Neal \cite[Section 1]{BlNe07}. Note, however, that the \nbd{C^*}subalgebra of $\cA$ generated by $F$ may but need not be isomorphic to the linking algebra. (Just take $F=\cA$ for $M_2(\cA)\not\cong\cA\ni\U$.) The statement in the lemma about the linking algebras is a statement about an off-diagonal corner of $M_2(\cA)$, not about $\cA$. It is important to never forget that, in these notes, we take the TRO structure of a Hilbert module from sitting in its linking algebra; the definition of TRO allows more general embeddings into $\cA$. The TRO-structure is independent of this choice (see Proposition \ref{TROuniprop} and Corollary \ref{isoisomcor}), but structures like the \nbd{C^*}algebra generated by $F$ or \it{positivity} (see the end of Problem \ref{Q3}) depend on the choice.
\eob

\bob \label{iliob}
(See, again, also Remark \ref{MERCrem}.) Matrices like linking algebras have been discussed in a more general context under the name of \it{generalized matrix algebras} in Skeide \cite{Ske00b} or \cite[Chapter 2]{Ske01}. The properties regarding ``\it{blockwise} subalgebra'' in the lemma are not automatic for subalgebras. For instance, if the linking algebra is unital (say, if $E$ is finite-dimensional), then the subalgebra $\C\U$ is not a \it{blockwise} or \it{matrix subalgebra}. A unital matrix subalgebra must contain at least the two units of each entry in the diagonal, separately. The more important is it to observe that an ideal of a matrix algebra is a matrix subalgebra, automatically. (Indeed, using approximate units for each diagonal entry from the left and from the right, we get projection maps $P_{i,j}$ onto each matrix entry.)
\eob

\proof[Proof of Theorem \ref{Ethm}.~]
Since \eqref{I1} and\eqref{I2} are equivalent by Proposition E, we only show equivalence of \eqref{I2} and\eqref{I3}. Again, for both directions (\eqref{I2}$\Rightarrow$\eqref{I3} also involving Corollary \ref{TI=TIcor}) we are left with straightforward verifications.\qed

\vspace{-1ex}
\bdefi
After the theorem, we say a \hl{closed ideal} in a Hilbert module $E$ is a subset $K$ fulfilling one (hence, all) of the conditions in the theorem.
\edefi

We add some \it{folklore} (see Remark  \ref{MERCrem}). By Observation \ref{iliob}, an ideal $\cI$ in the linking algebra of $E$ is a matrix subalgebra, that is, it has the form $\cI=\rtMatrix{I&K^*\\K&J}\subset\rtMatrix{\cB&E^*\\E&\sK(E)}$. If $\cI$ is closed, then necessarily: $I$ is a closed ideal in $\cB$; $K$ is a closed ternary ideal in $E$ and the ideal submodule associated with $I$ (clearly, $KI\subset K$, and an approximate unit for $\cI$ shows that $I$ acts nondegenerately); similarly, $K^*$ is the ideal submodule of the Hilbert \nbd{\sK(E)}module $E^*$ associated with the closed ideal $J$ of $\sK(E)$; moreover, since $E^*$ is full, (see the uniqueness statement in the proof of Proposition E), $J=\sK(E)_{K^*}=\sK(K)$. Summing up, every closed ideal $\cI$ in $\sK\rtMatrix{\cB\\E}$ gives rise to a unique closed ideal $I=P_{1,1}(\cI)$ in $\cB$ and every closed ideal $I$ in $\cB$ gives rise to a unique ideal $\cI$ in $\sK\rtMatrix{\cB\\E}$ such that $I=P_{1,1}(\cI)$. Requiring that $E$ is full, respectively, restricting to the reduced linking algebra, also $I$ is determined to be $\cB_K=\sK(K^*)$. We, thus, proved, the following:

\vspace{-1ex}
\bprop \label{extRCprop}
Let $E$ be a Hilbert \nbd{\cB}module.
\vspace{-1ex}
\begin{enumerate}
\item \label{RC1}
The formula $I=P_{1,1}(\cI)$ establishes a one-to-one correspondence between
\begin{enumerate}
\vspace{-1ex}
\item \label{RC1a}
closed ideals $I$ in $\cB$ and 

\item \label{RC1b}
closed ideals $\cI$ in $\sK\rtMatrix{\cB\\E}$,
\end{enumerate}
\vspace{-1ex}
also fulfilling $P_{2,1}(\cI)=\cls KI$ and $P_{2,2}(\cI)=\sK(\cls KI)$.

\lf
\item \label{RC2}
There are one-to-one correspondences between:
\begin{enumerate}
\vspace{-1ex}
\item \label{RC2a}
Closed ideals $\cI$ in the reduced linking algebra $\sK\rtMatrix{\cB_E\\E}$ of $E$.

\item \label{RC2b}
Closed ideals $I$ in $\cB_E$.

\item \label{RC2c}
Closed ideals $K$ in $E$.

\item \label{RC2d}
Closed ideals $J$ in $\sK(E)$.
\end{enumerate}
\vspace{-1ex}
The correspondences satisfy 
\baln{
P_{1,1}(\cI)
&
~=~
I
\;
~=~
\cls\AB{K,K},
\\
P_{2,1}(\cI)
&
~=~
K
~=~
\cls KI
~=~
\cls JK,
\\
P_{2,2}(\cI)
&
~=~
J
\,
~=~
\cls KK^*,
}\ealn
and are determined by (the appropriate subsets of) them (completed by $P_{1,2}(\cI)=K^*$).
\end{enumerate}
\eprop

\noindent
The symmetric situation in the second part is crucial for almost all that follows; not only because of the stated uniqueness properties (most of them desirable, if not indispensable), but also (and maybe even more importantly) by existence theorems for certain ``good'' maps in Section \ref{idealSEC}. The perfect symmetry between $E$ and $E^*$, $\cB_E=\sK(E^*)$ and $\sK(E)$ ($=(\sK(E))_{E^*}$) is part of \it{Morita equivalence} due to Rieffel \cite{Rie74,Rie74a}. In particular, some parts of Proposition \ref{extRCprop}\eqref{RC2} are known as \it{Rieffel correspondence}. However, as mentioned in the introduction, \it{Morita equivalence} and \it{Rieffel correspondence} come along in a different \it{clothing}. Since the explanation would distract us, we move a brief discussion together with other bibliographical hints to the following remark that ends this section. It may be skipped.

\brem \label{MERCrem}
If $E$ is a Hilbert \nbd{\cB}module, then $\sK(E)$ and $\cB_E$ are \it{Morita equivalent}. In fact, if we follow the definition  in Lance \cite[Chapter 7]{Lan95}, then two \nbd{C^*}algebras $\cA$ and $\cB$ are \hl{Morita equivalent} if there exists a full Hilbert \nbd{\cB}module such that $\sK(E)$ is isomorphic to $\cA$. In this point of view, the symmetry between $E$ and $E^*$, $\cB_E=\sK(E^*)$ and $\sK(E)$ explains why Morita equivalence is symmetric.

The original definition due to Rieffel \cite{Rie74,Rie74a} says $\cA$ and $\cB$ are (strongly) Morita equivalent, if there exists an \it{imprimitivity \nbd{\cA}\nbd{\cB}bimodule}.%
\footnote{
Rieffel's definition says \it{strongly Morita equivalent}, motivated by the fact that \it{strong Morita equivalence} is strictly stronger than what algebraic \it{Morita equivalence} asserts. It becomes more and more common to omit the word ``strong''.

Also, Rieffel's \it{imprimitivity bimodules} are defined for pre-\nbd{C^*}algebras (with bounded actions) and are not required to be complete. However, as Rieffel notes himself, the following results from Rieffel \cite[Section 3]{Rie79} are valid only for \nbd{C^*}algebras, and in this case the bimodule may be assumed to be complete. In definitions in more recent literature, \it{imprimitivity bimodules} occur for \nbd{C^*}algebras only, and are required to be complete.
}
Following Brown, Mingo, and Shen \cite{BMS94}, a \hl{Hilbert \nbd{\cA}\nbd{\cB}bimodule} is a Hilbert \nbd{\cB}module and a \bf{left} Hilbert \nbd{\cA}module such that $x\AB{y,z}=\AB{x,y}_\cA z$. (Equivalently, a Hilbert \nbd{\cA}\nbd{\cB}bimodule is an \nbd{\cA}\nbd{\cB}cor\-re\-spondence such that ``$\cA\supset\sK(E)$'' meaning that $\cA$ contains an ideal isomorphic to $\sK(E)$. In this case, $\AB{x,y}_\cA:=xy^*$.) An \hl{imprimitivity \nbd{\cA}\nbd{\cB}bimodule} is a Hilbert \nbd{\cA}\nbd{\cB}bimodule such that the underlying Hilbert module and left Hilbert module are full.%
\footnote{ \label{MEsymFN}
Optically, in Lance's definition of Morita equivalence, the symmetry between the algebra on the left, $\cA$, and the algebra on the right, $\cB$, is disturbed. Rieffel's original definition is more symmetric; but the price is that we have to deal with \bf{left} Hilbert modules. But, honestly, we think there is no real optical symmetry between the theory of Hilbert right modules (module map condition=associativity condition) and Hilbert left modules (module map conditions have to flip symbols in formulae); so we are not entirely happy calling this really symmetric. In Skeide \cite[Section 2]{Ske16} we propose a perfectly symmetric definition (using tensor products of correspondences), the way it is done in abstract algebra, and show its equivalence with (one of) the usual definitions.

Anyway, the story whether \it{left} or \it{right} is `nicer', starts, we think, in the very moment when we make our decision whether to write $f(x)$ or to write $(x)f$ for the value of a function $f$ at the point $x$ of its domain. (Apply this to the question whether the module-map-condition is a simple associativity or whether it has to flip around symbols.) This, actually, applies already to vector spaces, which, in a world where everybody writes $f(x)$, would be better off as right modules over their scalar field. And why would one prefer to write the linear functional $f_y$ on a Hilbert space induced by one of its elements $y$ in the form $f_y=\AB{\bullet,y}$ (flipping around the order of symbols in $f_y(x)=\AB{x,y}$), instead of $f_y=\AB{y,\bullet}$ (where in $f_y(x)=\AB{y,x}$ everything stays in order)?\vspace{1ex}
}
Of course, a Hilbert \nbd{\cA}\nbd{\cB}bimodule is also an \hl{\nbd{\cA}\nbd{\cB}correspondence} (that is, a Hilbert \nbd{\cB}module with a nondegenerate left action of $\cA$). By an \hl{\nbd{\cA}\nbd{\cB}Morita equivalence} we mean an imprimitivity \nbd{\cA}\nbd{\cB}bimodule viewed as \nbd{\cA}\nbd{\cB}correspondence.

What is known as the \hl{Rieffel correspondence} is the one-to-one correspondence between the ideals of $\cA$ and the ideals of $\cB$ provided $\cA$ and $\cB$ are Morita equivalent, and can be found in Rieffel \cite[Theorem 3.1]{Rie79}. Especially, when we use the Hilbert \nbd{\cB}module $E$ as an imprimitivity \nbd{\sK(E)}\nbd{\cB_E}bimodule, we recover equivalence of \ref{RC2b} and \ref{RC2d} in Proposition \ref{extRCprop}\eqref{RC2}. Using the isomorphism between $\cA$ and $\sK(E)$ as required in Lance's definition of Morita equivalence, we also recover the statement in Rieffel's generality. However, Rieffel's theorem states more: Firstly, the one-to-one correspondence is a lattice isomorphism under inclusion of ideals. (This follows simply by how this correspondence is constructed, be it in Rieffel \cite{Rie79}, be it in our proposition.) Secondly, and in our context more noteworthy, is a third equivalence stated in Rieffel's theorem but usually not mentioned in other works when referring to \it{Rieffel correspondence}. This is the equivalence of both lattices of ideals with the lattice of closed \nbd{\cA}\nbd{\cB}submodules of the imprimitivity bimodule. Closed \nbd{\cA}\nbd{\cB}submodules, when taking into account explicitly what full means for the left and the right action, correspond exactly to \it{ternary ideals} in the sense of \cite{BlNe07}; by Corollary \ref{TI=TIcor}, they correspond to our closed ternary ideals. This means, also equivalence with our \ref{RC2c} is (however, only after Corollary \ref{TI=TIcor}) included in \cite[Theorem 3.1]{Rie79}. Equivalence with \ref{RC2a} seems to be new (though as simple as the others, once properly formulated). And we are not aware to have seen Proposition  \ref{extRCprop}\eqref{RC1}, which is very important for our understanding of ideal submodules of not necessarily full Hilbert modules%
\footnote{
Note that in all our formulations in these notes, we carefully avoid to assume that the Hilbert modules are required to be full (to which most results in \cite{BaGu02,BaGu04,BaGu03} restrict). This is compensated by the fact that in our definitions we use the \bf{reduced} linking algebra, not the linking algebra. This reflects that in the \it{ternary} conditions we require, the actual ``size'' of $\cB\supset\cB_E$ does not play any role.
}%
, somewhere else.

\erem

\newpage
\section{Ideals and maps} \label{idealSEC}

In Definitions ISM and TI in  Section \ref{zero} and in Definition \ref{lidef}, we have defined three notions of (closed) ideal in a Hilbert module and we showed, in Theorem \ref{Ethm}, the (closed versions of the) three notions coincide. In this section we wish to justify the name ideal in Hilbert modules for this structure, by comparing it with what ideals are in \nbd{C^*}algebra theory. (The closely related notion of extension, we postpone to Section \ref{extSEC}.)

Closed ideals in a \nbd{C^*}algebra are precisely:
\begin{itemize}
\item
The subspaces that may be divided out, with a quotient in the same category.

\item
The kernels of the homomorphisms of the category.
\end{itemize}
The two are intimately related by:
\begin{itemize}
\item
If $I$ is a closed ideal in $\cB$, then the canonical map $b\mapsto b+I$ is a homomorphism with kernel $I$.

\item
If $\vp\colon\cB\rightarrow\cA$ is a homomorphism, then $\ker\vp$ is a closed ideal and $\vp(\cB)\cong\cB/\ker\vp$ via $\vp(b)\mapsto b+I$.
\end{itemize}
(The whole situation is also captured by the statement that we have a so-called \it{short exact sequence}
\beqn{ \tag{$*$} \label{*}
0
~\xrightarrow{~\text{can.}~}~
I
~\xrightarrow{~\text{can.}~}~
\cB
~\xrightarrow{~\text{can.}~}~
\cB/I
~\xrightarrow{~\text{can.}~}~
0,
}\eeqn
that is, $\cB$ is what is called an \it{extension} of $\cB/I$ by $I$. We discuss this later in Section \ref{extSEC}.)

So, our job is to find a quotienting procedure for closed ideals of Hilbert modules, to find the right sort of (homo)morphisms, and to give an interpretation in which these ideals are exactly the subspaces that can be divided out.

\lf
Actually, the search is not limited to the morphisms, but starts earlier with the question what is the class of objects. Are we considering a category of Hilbert \nbd{\cB}modules for fixed $\cB$ (for instance, like Hilbert spaces), or are we looking at a category of Hilbert modules over arbitrary \nbd{C^*}algebras? Are the modules required to be full or not? And is there a point of view where the question of fullness does not play a role? This has, obviously, to do with the question on which structure elements the three notions of ideal put emphasis.

Ideal submodules and the linking algebra put emphasis on $\cB$. For the former, everything is expressed and determined explicitly by choosing an ideal $I$ of $\cB$; the latter contains full information about $\cB$ in the \nbd{11}corner; different choices for $I$ giving the  same ideal submodule, remain visible only if we keep this information. (See Proposition \ref{extRCprop}\eqref{RC1}.) This suggests that ideal submodules, with full knowledge about the chosen ideal $I$, are best adapted to the situation where we wish to fix the \nbd{C^*}algebra $\cB$ and consider only Hilbert \nbd{\cB}modules. We see in a minute (Example \ref{qex}), that this is not really compatible with the idea of dividing out ideal submodules.

\lf
Let us reflect for a moment about quotients in general:

\begin{itemize}
\item
Of course, a Hilbert \nbd{\cB}module $E$ is a vector space and a closed ideal $K$ is a vector subspace. It is agreed that, when speaking about quotients, we start with the quotient vector space $E/K=\CB{x+K\colon x\in E}$ and, then, see what other structures it carries.

\item
$E$ is a normed (even a Banach) space and the quotient $E/K$ with a closed subspace $K$ is a normed (even a Banach) space for the quotient norm $\norm{x+K}:=\inf_{k\in K}\norm{x+k}$. Since we might define a norm on $E/K$ also in a different way, it is always a question if such an \it{ad hoc} norm coincides with the quotient norm.

\item
The quotient of a Hilbert space (so, $\cB=\C$) by a closed subspace with the quotient norm, is isomorphic to the orthogonal complement of the subspace with the norm calculated from the inherited inner product of the subspace. The same works for Hilbert modules as soon as the to-be-divided-out space is a complemented submodule.
\end{itemize}

\bprop \label{qprop}
Suppose the Hilbert \nbd{\cB}module $E$ decomposes as $E=F\oplus F^\perp$. Then for $q=q^*$, the projection onto $F^\perp$, the map $x+F\mapsto qx$ is an isomorphism $E/F\rightarrow F^\perp$.
\eprop

\proof
The map, once well defined, is clearly surjective. So, it remains to show that is is isometric. Since $x$ and $qx$ differ by an element in $F$, this follows from
\beqn{
\norm{qx+F}^2
~=~
\inf_{y\in F}\norm{qx+y}^2
~=~
\inf_{y\in F}\norm{\AB{qx,qx}+\AB{y,y}}
~=~
\norm{\AB{qx,qx}}
~=~
\norm{qx}^2,
}\eeqn
because $\norm{\AB{qx,qx}}\le\norm{\AB{qx,qx}+\AB{y,y}}\le\norm{\AB{qx,qx}}+\norm{\AB{y,y}}$.\qed

\lf\noindent
What can we say, if the submodule is not complemented? Well, in general:

\begin{itemize}
\item
Every quotient of a normed (Banach) module by a closed submodule is again a normed (Banach) module.

\item
By (a simple application of) Lance's \cite[Theorem 3.5]{Lan95}, there is at most one inner product on a Banach module that turns it into a Hilbert module. Actually, a careful analysis if the proof of \cite[Theorem 3.5]{Lan95} shows that this statement is true already for normed modules.
\end{itemize}

\noindent
For general closed submodules, we are not able to give an answer whether or not we are able to define an inner product on the quotient nor whether, if it exists, the norm induced by that inner product has to be the quotient norm.%
\footnote{
If we require the canonical map to be bounded (for instance, if the norm is the quotient norm), we could show that for non-complemented submodules, no such norm comes from an inner product, provided the following assertion was true: Every bounded right linear map on a Hilbert module that vanishes on a submodule $F$, also vanishes on $F^{\perp\perp}$. (Indeed, replacing $E$ with $F^{\perp\perp}$, we may assume that $F^\perp=\zero$. If the quotient is a Hilbert \nbd{\cB}module, then the canonical map is bounded, right linear, and vanishes on $F$, hence, on $E=F^{\perp\perp}$.) This is one more in a list of statements that could be proved if the assertion was true. (See the footnotes in Bhat and Skeide \cite[Section 2]{BhSk15}.) But despite some people saying it is true, we have not seen a working proof -- nor have we seen a counter example.
}
We now show that no norm on the quotient by a non-complemented ideal submodule is induced by a (\nbd{\cB}valued!) inner product.

\bex \label{qex}
Let $E$ be a (nonzero) full Hilbert \nbd{\cB}module and $K=\cls EI\ne E$ be the ideal submodule associated with a proper essential closed ideal $I$ in $\cB$ (so that $K\ne E$ and $K^\perp=\zero$). Suppose we had a \nbd{\cB}valued inner product on the quotient $E/K$. Choose a nonzero $x+K$ in the quotient, so, $b:=\sqrt{\AB{x+K,x+K}}\ne0$. Since $I$ is essential, we may choose $i\in I$ such that $bi\ne 0$. Then
\beqn{
\norm{(x+K)i}^2
~=~
\norm{i^*\AB{x+K,x+K}i}
~=~
\norm{(bi^*)(bi)}
~\ne~0,
}\eeqn
while $(x+K)i=xi+K=0\in E/K$. Consequently, there is no inner product on $E/K$ turning it into a (pre-)Hilbert \nbd{\cB}module.
\eex

It is noteworthy that (adding to both parts of Proposition \ref{extRCprop}) the correspondence between ideals $I$ in $\cB$ and ideal submodules $K=\cls EI$ in $E$ respects complements and direct sums, in the sense that that $K^\perp=\cls E(I^\perp)$. (Indeed, clearly $E(I^\perp)\subset (EI)^\perp$. And if $x\in(EI)^\perp$, then $\AB{x,y}i=0$ for all $y\in E$ and $i\in I$, so, $\AB{x,y}\in I^\perp$ for all $y$. Therefore, $\abs{x-xu_\lambda}^2\to0$ for every bounded approximate unit for $I^\perp$, so, $x\in\cls E(I^\perp)$.) With this, the statement of Example \ref{qex} extends to quotients with non-complemented ideals.

\bthm \label{qthm}
Let $K$ be a the closed ideal $K$ in the Hilbert \nbd{\cB}modules $E$. Then the quotient \nbd{\cB}module $E/K$ admits an inner product turning it isometrically into a Hilbert \nbd{\cB}module if and only if $K$ is complemented in $E$ (so that, necessarily, the norm is the quotient norm).
\ethm

\proof
We write $K=\cls EI$. If $K$ is complemented, then we are done by Proposition \ref{qprop}. So, for the converse, we assume $\norm{x+K}^2=\norm{\AB{x+K,x+K}}$ for some \nbd{\cB}valued inner product on $E/K$.

Clearly, the canonical map $x\mapsto x+K$ is injective on $K^\perp$. Since the inner product recovering $\norm{\bullet}\upharpoonright K^\perp$ is unique, the inner product of $E/K$ restricted to $K^\perp$ coincides with the original inner product of $K^\perp$. Therefore, $K^\perp$ is a complete, hence, closed ternary ideal also in $E/K$.

Since $(x+K)j=xj\in K^\perp$ for $j\in I^\perp$, we get $K^\perp=\cls (E/K)(I^\perp)$. Its complement in $E/K$ is the closure of $\ls(E/K)(I^{\perp\perp})$. The completion $\cls(E/K)(I^{\perp\perp})$ is a Hilbert \nbd{I^{\perp\perp}}module. Since $I$ is essential in $I^{\perp\perp}$, the closed ideal $\cls(E/K)I$ has zero complement. On the other hand, $\cls(E/K)I=\zero$. Therefore, $\cls(E/K)(I^{\perp\perp})=\zero$, hence, $E/K=K^\perp$, hence, $E=K\oplus K^\perp$.\qed

\lf
The key point in Example \ref{qex} for non-complemented $K$, is the contradiction that a \nbd{\cB}valued inner product had no choice but assigning  to certain  elements $x+K$ with $x\in K$ (so that $x+K=0$) a non-zero inner product $\AB{x+K,x+K}\ne0$. On the other hand, these non-zero elements are, actually, in $I$. If we divide out $I$ from $\cB$, then the preceding contradiction disappears. In fact, it is easy to show in general, that the \nbd{\cB/I}valued inner product $\AB{x+K,x'+K}=\AB{x,x'}+I$ is well-defined. Also, $E/K$ may be viewed as \nbd{\cB/I}module via $(x+K)(b+I)=xb+K$, and $E/K$ with this inner product is a Hilbert \nbd{\cB/I}module. (See Bakic and Guljas \cite{BaGu02} for details. See also Remark \ref{homrem}.) So, in order to quotient out an ideal submodule, we also have to take a quotient of the algebra $\cB$.

The quotient map $v\colon x\mapsto x+K$ shows us a sort of homomorphisms which have ideal submodules as kernels. In fact, denote by $\vp\colon b\mapsto b+I$ the quotient map for the quotient $\cB/I$. Then
\beq{ \label{phiiso}
\AB{vx,vx'}
~=~
\vp(\AB{x,x'}),
}\eeq
as observed in \cite{BaGu02} (see Remark \ref{homrem}), who called such a map a \it{\nbd{\vp}morphism}. Following Abbaspour and Skeide \cite{Ske06c,AbSk07}, we prefer to call a map $v$ from a Hilbert \nbd{\cB}module to a Hilbert \nbd{\cC}module fulfilling the preceding property for a homomorphism $\vp\colon\cA\rightarrow\cC$ a \hl{\nbd{\vp}isometry}; we say $v$ is a \hl{generalized isometry} if there exists a homomorphism $\vp$ turning it into a \nbd{\vp}isometry. (Note that a \nbd{\vp}isometry is norm-preserving if and only if $\vp$ is faithful on the range ideal of the \nbd{\cB}module.) We are sure that there will be earlier instances of \nbd{\vp}isometries and that most of their properties are \it{folklore}; in any case, most of their properties are easy exercises. (See, again, Remark \ref{homrem}.)

The canonical maps $K\rightarrow E\rightarrow E/K$ (that is, the analogues of those occurring in \eqref{*}) are generalized isometries. This depends, however, on the special choice of $K$ in the part $K\rightarrow E$ as a subspace of $E$ which \bf{is} a \nbd{\cB}module (so that $\vp=\id_\cB$ exists). We will see in Observation \ref{nonextob} that there are maps between Hilbert modules having ideal submodules as kernels, which are not generalized isometries.

\lf\noindent
\bf{Conclusion for \it{ideal submodules}.~}
\begin{itemize}
\item
The notion of ideal submodule, like the notion of generalized isometries as the corresponding morphisms, puts much emphasis on the algebra $\cB$. This might be felt as suggesting that we are looking at \nbd{\cB}modules for fixed $\cB$. On the contrary, the quotient by an ideal submodule (almost always) fails to admit an inner product with values in $\cB$.

\item
While the property of being an ideal submodule does not depend on the choice of the ideal $I$ that makes it explicit, the \nbd{C^*}algebra $\cB/I$ over which the quotient is a Hilbert module does depend on the choice of $I$. This is quite unfortunate.

\item
The derived notion of generalized isometry does not characterize all maps between Hil\-bert modules that have ideal submodules as kernels.
\end{itemize}
In conclusion, the notion ideal submodule is not the most fortunate choice to illustrate that ideals in Hilbert modules merit their name. (Actually, in order to give a completely asymmetric definition, according to Corollary \ref{TI=TIcor}, we might add to Theorem \ref{Ethm} the fourth equivalent property that $K$ be a closed \nbd{\cB}submodule of $E$ satisfying $\sK(E)K\subset K$. Also comparison with subbimodules of Hilbert bimodules as in Remark \ref{MERCrem} is worthwhile.)

\lf
Switching, now, our attention to ternary ideals, let us examine, first, in which sense closed ideals in \nbd{C^*}algebras are the \bf{unique} sort of subspaces that can be divided out. Already for vector spaces, we are used to the quotient set inheriting the vector space structure from the containing vector space: $(x+F)+(y+F)=(x+y)+F$ and $(x+F)\lambda=(x\lambda)+F$. (Actually, replacing the scalar $\lambda\in\C$ with $b\in\cB$, this also works for modules, resulting in the quotient of (normed) right \nbd{\cB}modules by their (closed) submodules.) For the quotient of algebras, we wish to define $(x+F)(y+F)=(xy)+F$. (This works if and only if $F$ is a two-sided ideal.) And if we have a \nbd{*}algebra, then we wish to define the involution $(x+F)^*=x^*+F$. (This works if and only if $F^*=F$.) And after ``the usual'' bit of work (for instance, like \cite[Theorem 3.1.4]{Mur90}): The quotient of a (pre-)\nbd{C^*}algebra by a closed \nbd{*}ideal, is a (pre-)\nbd{C^*}algebra under the quotient norm.

We just have convinced ourselves (Example \ref{qex}) that the quotient Banach module of a Hilbert \nbd{\cB}module by a closed submodule (like, for instance, an ideal submodule) is a Banach \nbd{\cB}module under the inherited operations, yes, but, in general, not a Hilbert \nbd{\cB}module. On the contrary, on the quotient $E/K$ of the Hilbert \nbd{\cB}module $E$, considered as TRO with the ternary product $(x,y,z)\mapsto(x\sds y\sds z):=x\AB{y,z}$, and a subspace $K$, taking also into account Corollary \ref{TI=TIcor}, we may (well-)define a ternary product via
\beqn{
(x+K\sds y+K\sds z+K)
~:=~
(x\sds y\sds z)+K
}\eeqn
if and only if $K$ is a ternary ideal.

Clearly, if $E/K$ with this ternary product \bf{is} a TRO, then, by definition of the ternary product on $E/K$, the canonical map is a \hl{ternary homomorphism}, namely, a linear map $v$ from a TRO $E$ to a TRO $F$ satisfying $v(x\sds y\sds z)=(vx\sds vy\sds vz)$. And: The kernel of every ternary homomorphism is, clearly, a ternary ideal.

So, after we show that $E/K$ is a TRO (possibly, in a unique way and under the quotient norm) and that ternary homomorphisms are bounded (so that their kernels are closed), we have established closed ternary ideals of TROs as \bf{the} analogue of closed ideals of \nbd{C^*}algebras, in the sense of the beginning of this section. The statements are true; however, their proofs fit much better the situation of linking ideals. Before summarizing the situation of ternary ideals and passing, then, to linking ideals (providing also the yet missing proofs), we state the following little uniqueness property. (A \hl{ternary isomorphism} is a bijective ternary homomorphism. After seeing, later on, that ternary isomorphisms are norm-preserving, the statement of the proposition also includes the metric structure of the TRO in question.)

Recall that, by definition, giving to a vector space $V$ a TRO structure means identifying it via a linear injection $u\colon V\rightarrow\cA$ with a TRO $uV$ in a \nbd{C^*}algebra $\cA$.

\bprop \label{TROuniprop}
Suppose the linear space $V$ has a ternary product $(\bullet\sds\bullet\sds\bullet)$ and suppose we have two linear injective maps $u_i$ from $V$ into \nbd{C^*}algebras $\cA_i$ such that $u_iV$ are TROs and $u_i(x\sds y\sds z)=(u_ix)(u_iy)^*(u_iz)$. Then $u\colon u_1x\mapsto u_2x$ is a ternary isomorphism. So, up to ternary isomorphism there is at most one TRO structure on $V$ reproducing the given ternary product.
\eprop

\proof
Straightforward verification.\qed

\lf\lf\noindent
\bf{Conclusion for \it{closed ternary ideals}.~}
\begin{itemize}
\item
The TRO-structure of a Hilbert module is perfect to capture its structure intrinsically, by which we mean without explicit reference to $\cB$.

\item
It is clear as crystal what the corresponding morphisms and closed ideals are, providing a perfect analogy with homomorphism and closed ideals in \nbd{C^*}algebras.

\item
But we left for later to prove a number of statements, because the proofs fit better into the next case.
\end{itemize}
In conclusion, the notion of closed ternary ideal resolves perfectly the problem we posed in the beginning of this section; but, the notion is not ``suggestive'' regarding how to, actually, prove it does.

\lf
We, now, come to the third point of view, linking ideals, where everything is defined in terms of the reduced linking algebra -- or, alternatively, of the linking algebra but requiring full modules. First of all, note that we could have defined linking ideals also by making reference to the linking algebra instead of the reduced linking algebra. As a condition, this is the same:

\bprop
For a subset $K$ of a Hilbert module $E$ the following conditions are equivalent:
\begin{itemize}
\item
$K$ is a closed linking ideal of $E$.

\item
$\rtMatrix{\sscls\AB{K,K}&K^*\\K&\sscls KK^*}$ is a closed ideal in the linking algebra of $E$.
\end{itemize}
\eprop
(We omit the simple proof that follows simply by observing that if $E$ is a Hilbert \nbd{\cB}module, then $\cB_E$ is an ideal in $\cB$ and, therefore, the reduced linking algebra of $E$ is an ideal in the linking algebra.) Also the proof of Theorem \ref{Ethm} does not depend on this choice. However, if we opted for the linking algebra instead of the reduced linking algebra, we would run into precisely the same difficulties that make ideal submodules to be not the most convenient choice: Too much emphasis on the \nbd{11}corner, when $E$ is non-full.

In the discussion leading to Proposition \ref{extRCprop} we have seen (taking also into account Observation \ref{iliob} and Lemma \ref{TsSlem}) that the closed ideals $\cI$ of the reduced linking algebra $\sK\rtMatrix{\cB_E\\E}=\rtMatrix{\cB_E&E^*\\E&\sK(E)}$ are uniquely parametrized by closed ternary ideals $K$ of $E$ as
\beqn{
\cI
~=~
\rsMatrix{\Scls\AB{K,K}&K^*\\K&\Scls KK^*}
~=~
\rsMatrix{\cB_K&K^*\\K&\sK(K)}
~=~
\sK\rsMatrix{\cB_K\\K},
}\eeqn
the reduced linking algebra of $K$. It is plain that also the quotient \nbd{C^*}algebra has the form
\beqn{
\sK\rtMatrix{\cB_E\\E}/\cI
~=~
\rsMatrix{\cB_E/\cB_K&E^*/K^*\\E/K&\sK(E)/\sK(K)},
}\eeqn
that $E/K$ is a ternary subspace of this quotient (identifying, now, the quotient as TRO), and that, since $E$ generates its reduced linking algebra, so does $E/K$ for the quotient (identifying the latter as the reduced linking algebra of $E/K$, $\sK\rtMatrix{\cB_E/\cB_K\\E/K}$. So, $E/K$ is a TRO and, since the norm on $\sK\rtMatrix{\cB_E\\E}/\cI$ is the quotient norm and since each corner projection is a contraction, also the norm on $E/K$ is the quotient norm. We summarize.

\bob
Under the correspondence \eqref{RC2c} $\Leftrightarrow$ \eqref{RC2a} in Proposition \ref{extRCprop} between closed ternary ideals $K$ in a Hilbert \nbd{\cB}module $E$ and closed ideals $\cI$ in the reduced linking algebra $\sK\rtMatrix{\cB_E\\E}$ of $E$:
\begin{itemize}
\item
The quotient $E/K$ is a TRO under the inherited ternary product and the quotient norm, and its linking algebra is $\sK\rtMatrix{\cB_E\\E}/\cI$.

\item
The proof, actually, goes by identifying $E/K$ as \nbd{21}corner in the latter quotient.
\end{itemize}
\eob

\lf\noindent
We now pass to homomorphisms, eventually, also providing the argument why, in the situation of Proposition \ref{TROuniprop}, also the norm of a TRO is determined by its algebraic structure (namely, the same way it is for a \nbd{C^*}algebra). While, by Observation \ref{iliob}, ideals in matrix algebras like the (reduced) linking algebra are matrix subalgebras, automatically, for homomorphisms the condition to act \hl{blockwise} (some say \it{corner preserving}) between matrix algebras (that is, sending the \nbd{ij}corner to the \nbd{ij}corner) is not.

\bex
Let $\cB$ be unital and consider the matrix subalgebra $\cC=\rtMatrix{\U\C&\\&\U\C}\subset\cA:=M_2(\cB)$. Let $U\in M_2$ be any unitary (scalar) matrix such that $U\bullet U^*$ does not leave invariant the diagonal. Then $\Phi:=(U\bullet U^*)\upharpoonright\cC$ is a homomorphism from the matrix algebra $\cC$ into the matrix algebra $\cA$ that is not onto a matrix subalgebra of $\cA$. Moreover, $U\bullet U^*$ itself is an automorphism of $\cA$ that is not blockwise.
\eex

Following our recipe to infer properties related to Hilbert modules from the analogue properties of their reduced linking algebras, for homomorphisms between reduced linking algebras we cannot think of any more natural condition but that it sends the \nbd{21}corner to the \nbd{21}corner. But this is already enough to force the homomorphism be blockwise.

\bprop \label{blockprop}
Let $\Phi\colon\rtMatrix{\cB_E&E^*\\E&\sK(E)}\rightarrow\rtMatrix{\cC_F&F^*\\F&\sK(F)}$ be a homomorphism between the reduced linking algebras of a Hilbert \nbd{\cB}module $E$ and a Hilbert \nbd{\cC}module $F$ such that $\Phi(E)\subset F$. Then $\Phi$ is blockwise.
\eprop

\proof
Clearly, $\Phi(E^*)=\Phi(E)^*\subset F^*$ and $\Phi(\sK(E))=\cls\Phi(E)\Phi(E)^*\subset\cls FF^*=\sK(F)$. Likewise, $\Phi(\cB_E)\subset\cC_F$.\qed

\lf
If $\Phi$ is blockwise, then we typically write $\Phi=\rsMatrix{\vp&u^*\\u&\phi}$. Since $\Phi$ is a \nbd{*}map, necessarily $u^*x^*=(ux)^*$.%
\footnote{
We ask the reader for forgiveness that we write homomorphisms between algebras with brackets ($\vp\colon b\mapsto\vp(b)$), while  we write the off-diagonal corners $u$ (where possible) without brackets (that is, rather as \it{operators} $u\colon x\mapsto ux$). Just imagine how ugly the proof of Theorem \ref{homothm} would look with more brackets ...
}

\vspace{-1ex}
\brem
More generally, by the same proof, a homomorphism between the linking algebras that preserves the \nbd{21}corner, preserves the reduced linking algebras (and acts blockwise between them). However, for the  linking algebras, the proposition may fail. Indeed, suppose $\cB_E\ne\cB$ is unital so that $\cB=\cB_E\oplus J$, and suppose that $\Phi(E)\ne F$ is complemented in $F$. Then we may extend a homomorphism between the reduced linking algebras by any homomorphism from $J$ to $\sK(\Phi(E)^\perp)\subset\sK\rtMatrix{\Phi(E)\\\Phi(E)^\perp}=\sK(F)$.
\erem

So, speaking about linking ideals, we now examine homomorphisms that are blockwise between reduced linking algebras. Most of the following theorem is \it{folklore}. (See Remark \ref{homrem}.) We prefer to give an independent and streamlined proof.

\bthm \label{homothm}
For a map $u$ from a full Hilbert \nbd{\cB}module $E$ to a Hilbert \nbd{\cC}module the following conditions are equivalent:
\vspace{-1ex}
\begin{enumerate}
\item \label{u1}
$u$ is a generalized isometry.

\item \label{u2}
$u$ is a ternary homomorphism.

\item \label{u3}
$u$ extends to a (blockwise, unique) homomorphism $\Phi=\rtMatrix{\vp&u^*\\u&\phi}\colon\rtMatrix{\cB&E^*\\E&\sK(E)}\rightarrow\rtMatrix{\cC&F^*\\F&\sK(F)}$.
\end{enumerate}
\vspace{-1ex}
Moreover, the equivalence between \eqref{u2} and \eqref{u3} remains true for arbitrary $E$, if we speak about homomorphisms $\Phi$ between the reduced linking algebras. In this case, $\Phi$ is surjective/injective if (and only if) $u$ is surjective/injective.
\ethm

\proof
Of course, \eqref{u3} implies \eqref{u1} (independently of whether $E$ is full or not, but we know from the proof of Proposition \ref{blockprop} that $\Phi$ is blockwise and unique, once $E$ is full), and \eqref{u1} implies \eqref{u2}. So, let us assume \eqref{u2} holds.

Let $a\in\sB^a(E)$ and define an operator on the pre-Hilbert \nbd{\cC}submodule $\ls(uE)(uE)^*F$ of $F$ by $(ux)(ux')^*y\mapsto(uax)(ux')^*y$. It is easy to show that the operator corresponding to $a^*$ is a formal adjoint on the generating set $(uE)(uE)^*F$. Therefore, this well-defines a representation of the \nbd{C^*}algebra $\sB^a(E)$ by adjointable operators on $\ls(uE)(uE)^*F$. Since the unital \nbd{C^*}algebra $\sB^a(E)$ is spanned by its unitaries, the representation operators are bounded, so that the representation extends further to a representation by bounded operators on $F_E=\cls(uE)(uE)^*F\subset F$ mapping, clearly, $\sK(E)$ into $\sK(F_E)\subset\sK(F)$. We denote this map by $\phi$. Doing the same for the ternary homomorphism $u^*\colon x^*\mapsto(ux)^*$ from the Hilbert \nbd{\sK(E)}module $E^*$ to the Hilbert \nbd{\sK(F)}module $F^*$, we obtain a homomorphism $\vp$ from $\cB=\cB_E=\sK(E^*)$ into $\sK(F^*)=\cC_F\subset\cC$. It is routine to show that the matrix $\Phi$ is a homomorphism.

$\sK\rtMatrix{\cB_E\\E}$ and $\sK\rtMatrix{\cC_F\\F}$ are generated as \nbd{C^*}algebras by $E$ and $F$, respectively. So, $\Phi$ is surjective if and only if $u$ is. The homomorphism $\phi\colon\sK(E)\rightarrow\sK(F)$ has been constructed as the unique map $\sK(E)\rightarrow\sK(\cls (uE)(uE)^*F)\subset\sK(F)$ acting as $a\colon(ux)\AB{ux,y}\mapsto(u(ax))\AB{ux,y}$. Clearly, $\phi(a)(ux)\AB{u(ax),u(ax)}=(u(ax))\AB{u(ax),u(ax)}$ is $0$ if and only if $u(ax)$ is $0$, so $\phi(a)$ is $0$ if and only if $u(ax)$ is $0$ for all $x\in E$. A similar statement follows for $u^*$ and $\vp$. So, $\Phi$ (in Theorem \ref{homothm}) is injective if and only if $u$ is.\qed

\lf
And, now, here is the  missing ``normed'' supplement to Proposition \ref{TROuniprop}:

\bcor \label{isoisomcor}
Ternary isomorphisms are isometric.
\ecor

\proof
Let $u\colon E\rightarrow F$ and $v\colon F\rightarrow E$ be an inverse pair of ternary homomorphisms. Then their extensions $\Phi$ and $\Psi$ to homomorphisms between the reduced linking algebras are inverse to each other, too. Algebraic isomorphisms between \nbd{C^*}algebras are isometric, hence, \it{a fortiori}, so are their restrictions $u$ and $v$.\qed

\bob \label{nonextob}
 If $E$ is not full, then \eqref{u2} or \eqref{u3} do not imply \eqref{u1}. For instance, if $I$ is a proper essential ideal in $\cB$ then the identity on $I$ does not extend to a homomorphism $\cB\rightarrow I$. The ternary homomorphism $u=\id_I$ from the Hilbert \nbd{\cB}module $I$ to the Hilbert \nbd{I}module $I$ does not extend to a homomorphism between the linking algebras.
\eob

\bob
Generalized homomorphisms $E\rightarrow F$ correspond one-to-one with blockwise homomorphisms between the linking algebras, in a way not dissimilar to how ideals in $\cB$ correspond to ideals in the linking algebra of $E$ in Proposition \ref{extRCprop}\eqref{RC1}. However, while the \nbd{11}corner $I$ alone determines the ideal in the linking algebra, for homomorphisms between the linking algebras, we need to know the \nbd{11}corner $\vp$ and the \nbd{21}corner $u$ of $\Phi$. Given only $\vp$, there need not even exist a \nbd{\vp}isometry $u$ for given $E$ and $F$, hence, no $\Phi$; see \cite{Ske06c} for more details. Also for the reduced linking algebras, Conditions \eqref{RC2b} and \eqref{RC2c} in Proposition \ref{extRCprop}\eqref{RC1} that refer to the diagonal corners in Proposition \ref{extRCprop}\eqref{RC1}, have no counterpart for homomorphisms. But we do have a one-to-one correspondence between ternary homomorphisms and blockwise homomorphisms between the reduced linking algebras.
\eob

\lf\lf\noindent
\bf{Conclusion for \it{closed linking ideals}.~}
\begin{itemize}
\item
Linking ideals with maps that extend as blockwise homomorphisms between the reduced linking algebras perfectly resolve the problem of dividable submodules and homomorphisms posed in the beginning of this section.

\item
However, it turns out that when we want to characterize the homomorphisms without actually having to construct their extensions, then we end up with exactly the ternary homomorphisms.
\end{itemize}

\lf\noindent
\bf{Summary of conclusions.~}
\begin{itemize}
\item
The notions that perfectly resolve the quotient/homomorphism-problem in an intrinsic way, that is, requiring not more than looking at the modules and maps between them, are the notions of closed ternary ideal and of ternary homomorphism.

\item
The notions of ideal submodule and generalized isometry explicitly involve other corners of the linking algebra and, at least for non-full $E$, there may be unfortunate choices that ruin a satisfactory solution of our problem.

\item
While the ternary notions resolve our problem intrinsically, when we wish to show they do, we had best run through linking notions (closed linking ideal and blockwise homomorphisms), based on our theorems that establish their equivalence.
\end{itemize}

\brem
In many situations it has proved fruitful to introduce properties of Hilbert modules and of maps between Hilbert modules by passing to the (reduced) linking algebras and examining how ``good'' blockwise properties of the linking algebras are reflected by those of their \nbd{21}corners. (Only for instance: \it{Von Neumann modules} are Hilbert modules over a von Neumann algebra whose linking algebra is a von Neumann algebra and maps between them are \it{normal} if they possess a normal (in particular, positive!) blockwise extension; see \cite{Ske00b}. Only the definition of \it{dynamical systems on Hilbert modules} as one parameter groups of ternary automorphisms allowed to provide a neat characterization of their generators as \it{ternary derivations}, while their characterization as \it{generalized derivation} fails due to domain problems; see \cite{AbSk07}. Only studying the extension of so-called \nbd{\tau}maps (a \nbd{\tau}isometry but with $\tau$ a CP-map; see Question \ref{Q2}) to a blockwise map between the linking algebras allows to give an intrinsic characterization without reference to $\tau$; see \cite{SkSu14}.)

In the following section we add as an example the notion of \it{extensions} for Hilbert modules. Extensions have been introduced by Bakic and Guljas \cite{BaGu04} (examined in their forthcoming papers and also by Kolarec), based on the notion from \cite{BaGu02}, which correspond to the first of the three points of view we discussed until here. Without any ambition to be complete, we recover some of the results on extensions of Hilbert modules as corollaries of the known results on extensions of \nbd{C^*}algebras.
\erem

\brem \label{homrem}
It is difficult to sort out what in the present section was new or where exactly it has occurred before.

When speaking about quotients, surely (under the same ``translation'' as explained in Remark \ref{MERCrem}) one can trace back a good number of statements to Rieffel \cite[Section 3]{Rie79} (in particular, \cite[Corollary 3.2 and Proposition 3.2]{Rie79}).

As for the maps that occur as homomorphisms, also ternary homomorphisms necessarily have been around, explicitly or implicitly, since ternary rings (algebraic, \nbd{C^*}rings, or TRO) existed. (See  Remark \ref{MERCrem}.) With reference to Hilbert modules they might have occurred first in Bakic and Guljas \cite{BaGu02}. Most of Theorem \ref{homothm} occurs earlier. Equivalence of  \eqref{u1} and \eqref{u2} appears, for instance, in Abbaspour and Skeide \cite[Theorem 2.1]{AbSk07}.%
\footnote{
We mentioned already in Skeide and Sumesh \cite[Footnote 2]{SkSu14} that, unlike what is stated in \cite[Theorem 2.1]{AbSk07}, the hypothesis that a ternary homomorphism has to be linear cannot be dropped.
}
And results allowing to conclude that \eqref{u1} implies \eqref{u3} (the other direction being obvious) have been around in \cite{BaGu02,Ske06c,AbSk07}. In a sense, the problem in \eqref{u1} $\Rightarrow$ \eqref{u3} (namely, constructing the \nbd{22}corner of the extension $\Phi$ of a \nbd{\phi}isometry $u$ to the reduced linking algebra), is very much similar to constructing the extension of a ternary homomorphism $u\colon E\rightarrow F$ to the \nbd{11}corners as \eqref{u2} $\Rightarrow$ \eqref{u1}. (The proofs are quite different, though. The proof of \cite[Theorem 2.15]{BaGu02} is based on conjugation with a \nbd{\vp}unitary and reduction of the non-injective case to this by a quotient, and needs preparation. The proof we present here is self-contained. It is, essentially, the same proof as in \cite{AbSk07} but exchanging $\cB$ with $\sB^a(E)$ and $E^*$ with $E$.) In any case, probably none of these references is primary. 
Corollary \ref{isoisomcor} is Blecher and Le Merdy \cite[Lemma 8.3.2(2)]{BlLMe04}. (Of course, also \cite[Lemma 8.3.2(1)]{BlLMe04}, stating that ternary homomorphisms are contractions, is a simple consequence of Theorem \ref{homothm}, and that the statements extend to completely contractive/isometric is not an issue.) They attribute the lemma to the works of Harris and Kaup (which we did not check; \cite[Section 8.7]{BlLMe04}). The proof in \cite{BlLMe04} is by spectral calculus and direct. Ours here, appeals to the uniqueness of a \nbd{C^*}norm (also a ``very spectral'' property) on the linking algebra. An interesting related result by Solel \cite{Sol01} characterizes the (isometric) Banach space isomorphisms between Hilbert modules as exactly those maps that extend to a blockwise Banach space isomorphism between the reduced linking algebras that is the sum of a homomorphism and an anti-homomorphism.

It is clear that using the linking algebra to prove results about Hilbert modules is a very well known technique. (Though, surely, there are different ways to put this idea into practise to prove the same result.) What, we think, is novel -- if not to these notes, then surely to our earlier work \cite{Ske00b,AbSk07,SkSu14} -- is the idea to \it{motivate} notions for Hilbert modules by applying well known notions for \nbd{C^*}algebras to their linking algebras, and to use this to derive the properties of the former as corollaries of the latter. And these notes add further instances for this concept.

Last but not least, we emphasize once more that our Hilbert modules take their TRO structure from sitting as \nbd{21}corner in their reduced linking algebras, not from sitting as a general ternary subspace of a general \nbd{C^*}algebra $\cA$. (See, again, Observation \ref{TROob} for the difference.) All our extension theorems are for the reduced linking algebras, not for the \nbd{C^*}algebras generated by a TRO sitting in $\cA$.
\erem

\newpage
\section{Extensions} \label{extSEC}

An \it{extension} of \nbd{C^*}algebras is essentially a short exact sequence
\beqn{ \tag{$**$} \label{**}
0
~\xrightarrow{~0~}~
\cA
~\xrightarrow{~\psi~}~
\cB
~\xrightarrow{~\vp~}~
\cC
~\xrightarrow{~0~}~
0.
}\eeqn
A sequence of maps is \hl{exact} if the range of each map coincides with the kernel of the next. An exact sequence of three, augmented by the $0$s in the diagram, is \hl{short}. Since we are speaking of \nbd{C^*}algebras, the maps are homomorphisms. $0$ as \nbd{C^*}algebra is short for $\zero$. And the only (linear) maps that have $\zero$ as domain or as codomain, are the zero-maps $0$ that map everything to $0$. (Usually, we will not write the $0$s on arrows.) The situation of this short exact sequence is (usually) referred to by saying, $\cB$ is an \hl{extension} of $\cC$ by $\cA$.

Since, by exactness, the kernel of $\psi$ is $\zero$, the homomorphism $\psi$ is an embedding and $\cA$ is isomorphic to the subalgebra $\psi(\cA)$ of $\cB$. Since, by exactness, $\psi(\cA)=\ker\vp$, $\psi(\cA)$ is an ideal in $\cB$. Since, by exactness, $\vp$ is surjective, we see that $c=\vp(b)\mapsto b+\psi(\cA)$ establishes an isomorphism $\cC\cong\cB/\psi(\cA)$. So, \eqref{*} does, indeed, capture (up to isomorphism) the general situation of extension. However, while \eqref{*} puts emphasis on the algebra called extension, $\cB$, in extension theory one would rather fix $\cA$ and $\cC$ and analyze the class of all extensions of $\cC$ by $\cA$. (Also for this reason it is reasonable to use for $\cA$ a letter for generic \nbd{C^*}algebras, not a letter for ideals.)

\bdefi \label{TEdefi}
An \hl{extension} (or \hl{ternary extension}) of Hilbert modules is a short exact sequence
\beqn{ \tag{$*\!\!\!\;*\!\!\!\;*$} \label{***}
0
~\xrightarrow{~~~}~
G_\cA
~\xrightarrow{~v~}~
E_\cB
~\xrightarrow{~u~}~
F_\cC
~\xrightarrow{~~~}~
0
}\eeqn
of ternary homomorphisms. We say $E$ is an extension of $F$ by $G$.
\edefi
(The notation $E_\cB$ means $E$ is a Hilbert \nbd{\cB}module. We choose a ``lettering'', that is most compatible with both the notation from the preceding sections and the above short exact sequence of \nbd{C^*}algebras.)

The modules are not required to be full; fullness does not play a role in our definitions. However, we have the following obvious consequence of Proposition E in Section \ref{zero} and of Theorem \ref{homothm}:

\bprop
For full Hilbert modules the definition of extension is equivalent to the definition in \cite{BaGu04}. (For non full modules it is undefined in \cite{BaGu04}.)
\eprop

\bthm \label{extthm}
Extensions
\vspace{-1ex}
\beqn{
0
~\xrightarrow{~~~}~
G
~\xrightarrow{~v~}~
E
~\xrightarrow{~u~}~
F 
~\xrightarrow{~~~}~
0
}\eeqn
of Hilbert modules are in one-to-one correspondence with \hl{blockwise} extensions (that is, all homomorphisms are blockwise)
\beqn{
0
~\xrightarrow{~~~}~
\sK\rsMatrix{\cA_G\\G}
~\xrightarrow{~\Psi~}~
\sK\rsMatrix{\cB_E\\E}
~\xrightarrow{~\Phi~}~
\sK\rsMatrix{\cC_F\\F}
~\xrightarrow{~~~}~
0
}\eeqn
of the corresponding reduced linking algebras, via the requirement that the (co)restriction of $\Phi$ and $\Psi$ to maps between the \nbd{21}corners are $u$ and $v$, respectively.
\ethm

\proof
This is the one-to-one correspondence \eqref{u2}$\Leftrightarrow$\eqref{u3} in Theorem \ref{homothm} (including injectivity of $\Psi$ and surjectivity of $\Phi$) plus the statement that $\Psi$ is onto the kernel of $\Phi$ if (and only if) $v$ is onto the kernel of $u$.

For the latter, observe that $\sK\rtMatrix{\cB_{vG}\\vG}$ is an (automatically blockwise) ideal in $\sK\rtMatrix{\cB_E\\E}$ (because $vG$ is a ternary ideal in $E$), and that $\Phi$ is $0$ on that ideal (because $u\upharpoonright vG$ is). So, $\Phi$ gives rise to a (blockwise) homomorphism $\sK\rtMatrix{\cB_E\\E}/\sK\rtMatrix{\cB_{vG}\\vG}\rightarrow\sK\rtMatrix{\cC_F\\F}$. This homomorphism is injective, because the quotient map $E/vG\rightarrow F$ induced by $u$ is injective. In other words, $\Psi$ is onto $\ker\Phi$.\qed

\bcor
By (co)restricting the corresponding blockwise extension of an extension $E$ of $F$ by $G$ to the two diagonal corners, this extension gives rise to an associated extension $\cB_E$ of $\cC_F$ by $\cA_G$ and an associated extension $\sK(E)$ of $\sK(F)$ by $\sK(G)$.
\ecor

We follow, as we said without any ambition of being complete, through a (very) few basic properties of extensions of \nbd{C^*}algebras and their module analogues. This should be well readable together with Blackadar \cite[Section II.8.4]{Bla06}. It is understood that when we say blockwise extension we have Hilbert modules around as in Theorem \ref{extthm}, and we are referring to the line with the corresponding reduced linking algebras.

\bemp[Split extensions.~]
The extension \eqref{**} of \nbd{C^*}algebras is \hl{split} if there exists a homomorphism $s\colon\cC\rightarrow\cB$ (called a \hl{splitting}) that is a right inverse of $\vp$. (So, $s\circ\vp$ is a conditional expectation onto $s(\cC)\cong\cC$; see Question \ref{Q1}.) This definition makes sense for extensions of Hilbert modules, if we replace the homomorphism $s$ for \eqref{**} with a ternary homomorphism $s\colon F\rightarrow E$ for \eqref{***}. It is an immediate consequence of Theorem \ref{extthm} that this requirement is equivalent to requiring that the corresponding blockwise extension is split via a blockwise splitting (being the unique blockwise homomorphism corresponding to $s$ via Theorem \ref{homothm}).

Also the \hl{super-trivial}%
\footnote{ \label{stFN}
There is an equivalence relation among extensions that explicitly quotients out the split extensions. (See \cite[Definition II.8.4.18]{Bla06}.) Under this quotient, of course, all split extensions coincide with the super-trivial one. This is why the split extensions are referred to as the \hl{trivial} ones. But, here we mean a special representative among the trivial extensions, and call it the super-trivial one.)
}
extension $\cC\oplus\cA$ of $\cC$ by $\cA$ (with $\vp$ being the canonical projection and $\psi$ being the canonical embedding) generalizes. To see that, recall that the \hl{external} direct sum of a Hilbert \nbd{\cC}module $F$ and a Hilbert \nbd{\cA}module $G$ is $F\oplus G$ with its natural Hilbert \nbd{\cC\oplus\cA}module structure (so that also $\sK(F\oplus G)=\sK(F)\oplus\sK(G)$). It follows that $F\oplus G$ is an extension (obviously split) of $F$ by $G$.
\eemp

\bemp[Busby invariant.~]
The \hl{Busby invariant} of an extension $\cB$ of $\cC$ by $\cA$ as in \eqref{**} is the canonical homomorphism\vspace{-1.5ex}
\beqn{
\tau
\colon
\cC
~\xrightarrow{~\cong~}~
\cB/\psi(\cA)
~\longrightarrow~
M(\cA)/\cA
~=:~
Q(\cA),\vspace{-.5ex}
}\eeqn
where the second arrow is induced from the canonical homomorphism from $\cB$ into the multiplier algebra of its ideal $\psi(\cA)\cong\cA$. (Recall: If $\vt\colon\cD\rightarrow\mathcal{E}$ is a homomorphism and if $I$ is an ideal in $\cD$ such that $\vt(I)$ is an ideal in $\mathcal{E}$, then $d+I\mapsto\vt(d)+\vt(I)$ well-defines a homomorphism $\cD/I\rightarrow\mathcal{E}/\vt(I)$. Here, $I$ is $\psi(\cA)$ and the image of $I$ in $M(\cA)$ under $\vt$ is $\cA$.)

The only thing it takes to make work the definition of Busby invariant for an extension \eqref{***} of Hilbert modules as the \nbd{21}corner of the Busby invariant of the corresponding blockwise extension of the reduced linking algebras, is to get hold of the blockwise structure of the multiplier algebra $M\family{\sK\rtMatrix{\cA_G\\G}}$ and the \hl{corona} $Q\family{\sK\rtMatrix{\cA_G\\G}}$ of the reduced linking algebra $\sK\rtMatrix{\cA_G\\G}$. So, starting from an extension \eqref{***}, by Theorem \ref{extthm} we get the (unique) blockwise extension
\beqn{ 
0
~\xrightarrow{~~~}~
\sK\rsMatrix{\cA_G\\G}
~\xrightarrow{~\Psi~}~
\sK\rsMatrix{\cB_E\\E}
~\xrightarrow{~\Phi~}~
\sK\rsMatrix{\cC_F\\F}
~\xrightarrow{~~~}~
0
}\eeqn
(such that the (co)restrictions to the \nbd{21}corner of $\Phi$ and $\Psi$ are $u$ and $v$, respectively). For calculating its Busby invariant $\sT\colon\sK\rtMatrix{\cC_F\\F}\rightarrow Q\family{\sK\rtMatrix{\cA_G\\G}}=M\family{\sK\rtMatrix{\cA_G\\G}}/\sK\rtMatrix{\cA_G\\G}$, recall that for any Hilbert module $E$ we have $M(\sK(E))=\sB^a(E)$. (Kasparov \cite{Kas80}; see Skeide \cite[Corollary 1.7.14]{Ske01} for a simple proof.) So,
\beqn{
M\family{\sK\rsMatrix{\cA_G\\G}}
~=~
\sB^a\rsMatrix{\cA_G\\G}
~=~
\rsMatrix{M(\cA_G)&\sB^a(G,\cA_G)\\\sB^a(\cA_G,G)&\sB^a(G)}
~\supset~
\rsMatrix{\cA_G&G^*\\G&\sK(G)}
~=~
\sK\rsMatrix{\cA_G\\G}.
}\eeqn
It follows that
\beqn{
Q\family{\sK\rsMatrix{\cA_G\\G}}
~=~
\rsMatrix{M(\cA_G)/\cA_G&\sB^a(G,\cA_G)/G^*\\\sB^a(\cA_G,G)/G&\sB^a(G)/\sK(G)}
~=~
\rsMatrix{Q(\cA_G)&Q(G)^*\\Q(G)&Q(\sK(G))},
}\eeqn
where we defined $Q(G):=\sB^a(\cA_G,G)/G$. Note that this is really the quotient of the Hilbert \nbd{M(\cA_G)}module $M(G):=\sB^a(\cA_G,G)$ and its ternary ideal $G=\sK(\cA_G,G)$.

\bobn
If $\cA_G$ is nonunital, then $M(G)$ need not be full but only what we call \hl{strictly full}. The most drastic way to see this, is if $\sK(G)$ is unital (that is, if $G$ is algebraically finitely generated; for instance, if $G$ is the full Hilbert \nbd{\sK(H)}module $H^*$ for some infinite-dimensional Hilbert space $H$, so that $G$ is algebraically singly generated by any of its nonzero elements). In that case, $M(G)=\sB^a(\cA_G,G)=G$. (The easiest way to see this is, to consider the Hilbert module $G^*$ over the unital \nbd{C^*}algebra $\sK(G)$.) Then, $M(\cA_G)_{M(G)}=\cls\AB{M(G),M(G)}=\cls\AB{G,G}=\cA_G\subsetneq M(\cA_G)$. So, looking at $M\family{\sK\rtMatrix{\cA_G\\G}}$, we are definitely leaving the situation where our matrix algebras are reduced linking algebras. Even worse, in the same situation we have $Q\family{\sK\rtMatrix{\cA_G\\G}}=\rtMatrix{Q(\cA_G)&0\\0&0}$. A similar thing happens, if $G$ is full over a unital algebra with nonunital $\sK(G)$ (for instance $G=H$); just that now the zero-place and nonzero-place in the diagonal switch. Taking the external direct sum of these two cases, we get an example where both $\cA_G$ and $\sK(G)$ are nonunital and where $Q(G)=\zero$ so that $Q\family{\sK\rtMatrix{\cA_G\\G}}$ is diagonal with both places in the diagonal nonzero.
\eobn

In view of this observation, it appears appropriate to introduce the \hl{reduced multiplier linking algebra} $M^{red}(G)$ of a Hilbert \nbd{\cA}module $G$ as the \nbd{C^*}subalgebra of $M\family{\sK\rtMatrix{\cA_G\\G}}$ generated by $M(G)$. Of course, $M ^{red}(G)\supset\sK\rtMatrix{\cA_G\\G}$ sits as an ideal in $M\family{\sK\rtMatrix{\cA_G\\G}}$ and the \hl{reduced corona linking algebra} $Q^{red}(G):=M ^{red}(G)/\sK\rtMatrix{\cA_G\\G}\subset Q(\sK\rtMatrix{\cA_G\\G})$ sits as an ideal in $Q(\sK\rtMatrix{\cA_G\\G})$ and coincides with the \nbd{C^*}subalgebra generated by $Q(G)$. Since both $M^{red}(G)$ and $Q^{red}(G)$ are generated by their \nbd{21}corners $M(G)$ and $Q(G)$, respectively, $M^{red}(G)$ is the reduced linking algebra of $M(G)$ and $Q^{red}(G)$ is the reduced linking algebra of $Q(G)$. Clearly, a (blockwise) homomorphism $\sT\colon\sK\rtMatrix{\cC_F\\F}\rightarrow Q^{red}(G)$ may be considered a (blockwise) homomorphism $\sT\colon\sK\rtMatrix{\cC_F\\F}\rightarrow Q\family{\sK\rtMatrix{\cA_G\\G}}$.

Summarizing: We have a one-to-one correspondence between between
\begin{itemize}
\item[(i)]
ternary homomorphisms $T\colon F\rightarrow Q(G)$ and

\item[(ii)]
blockwise homomorphisms $\sT=\rtMatrix{\tau&T^*\\T&\vt}\colon\sK\rtMatrix{\cC_F\\F}\rightarrow Q^{red}(G)\subset Q\family{\sK\rtMatrix{\cA_G\\G}}$
\end{itemize}
such that the \nbd{21}corner of $\sT$ is $T$. If $\sT$, considered as map into $Q\family{\sK\rtMatrix{\cA_G\\G}}$, is the Busby invariant of the blockwise extension as in Theorem \ref{extthm}, then we say $T$ is the \hl{Busby invariant} of the corresponding extension of Hilbert modules.

\bobn
\begin{enumerate}
\item \label{Bob1}
The \it{multiplier space} $M(G)$ of $G$ shares (expressed in terms of ternary ideals and morphisms) the same universal property as the multiplier algebra of a \nbd{C^*}algebra. (Also this may be deduced by applying the universal property of the multiplier algebra of the reduced linking algebra of $G$.) Therefore, the homomorphism $T\colon F\rightarrow E/vG\rightarrow Q(G)$ maybe defined directly from an extension \eqref{***} without passing through the blockwise extension, and the result, obviously, coincides. (This direct definition by analogy with the case of \nbd{C^*}algebras, is what has been done in \cite{BaGu03} in terms of generalized isometries and full modules (which, we know, is equivalent to the procedure in terms of ternary homomorphisms we just discussed).) Our definition above recovers the Busby invariant from the Busby invariant for blockwise extensions of reduced linking algebras.

\item \label{Bob2}
Every homomorphism $\cC\rightarrow Q(\cA)$ is the Busby invariant of an extension $\cB$ of $\cC$ by $\cA$; thus, the Busby invariant induces an equivalence relation among extension (usually, referred to as \it{strong isomorphism}).%
\footnote{
This relation is stronger than the one we mentioned in Footnote \ref{stFN}. It does not identify split extensions with the super-trivial one.
}
Adding the one-to-one correspondence between extension of Hilbert modules and blockwise extensions of the reduced linking algebras, we get the same statement for the Busby invariants of extensions of Hilbert modules.

The proof for \nbd{C^*}algebras is done by proving existence of so-called \it{pullbacks}. Again, one may (as in \cite{BaGu03}) imitate word by word the \nbd{C^*}case. Or one may note that, for blockwise homomorphisms, the construction of a pullback for reduced linking algebras, leads to blockwise homomorphisms and, therefore, restricts to a pullback construction for the \nbd{21}corners. We do not deepen this further, but leave this as an exercise. And we spare the discussion of equivalences of extensions. (See \cite{BaGu03}.) We only keep the message: Every ternary homomorphism $F\rightarrow Q(G)$ is the Busby invariant of an extension $E$ of $F$ by $G$.

\item \label{Bob3}
We briefly address the the situation when the Busby invariant $T$ of an extension of Hilbert modules \eqref{***} is $0$, and use this to be more explicit about how the \it{canonical homomorphism} from $E$ into $M(G)$ looks like. By \eqref{***} being a short exact sequence, it follows that $vG$ sits as a closed ternary ideal in $E$. So let us assume for the time being, that $G$ actually is a closed ternary ideal in $E$ and that $v$ is the canonical identification of $G\subset E$. Recall that then $\cA_G$ is a closed ideal in $\cB_E$, and $G=\cls E\cA_G$. The multiplier space of $G$ is $M(G)=\sB^a(\cA_G,G)$ and every element $x$ in $E$ gives rise to a map $\sB^a(\cA_G,G)\ni x\colon a\rightarrow xa$. The Busby invariant is $0$ if and only if the map $E\rightarrow\sB^a(\cA_G,G)$ is into $G\subset\sB^a(\cA_G,G)$. In this case, for each $x\in E$ there is a unique (since $G$ is full) $g_x\in G$ such that $xa=g_xa$ for all $a\in\cA_G$. By setting $px:=g_x$, we define a map $p\colon E\rightarrow  G$. Since $\AB{E,G}\subset\cA_G$, by choosing an approximate unit $\bfam{u_\lambda}_{\lambda\in\Lambda}$ for $\cA_G$, we get
\beqn{
\AB{x,px'}
~=~
\lim_\lambda\AB{xu_\lambda,px'}
~=~
\lim_\lambda\AB{pxu_\lambda,px'}
~=~
\AB{px,px'}.
}\eeqn
It is an intriguing exercise to show that this implies that $p$ is adjointable (hence, linear) and fulfills $p^*p=p$; in other words, $p$ is a projection. It follows that $G$ is a direct summand in $E$, so $E=G\oplus G^\perp$. As explained ahead of Theorem \ref{qthm}, this means that also $G^\perp=\cls E(\cA_G^\perp)$ is an ideal and that $\cB=\cA_G\oplus\cA_G^\perp$. By Proposition \ref{qprop}, $E/G\cong G^\perp$, so $F\cong G^\perp$.
Returning to the general situation, we have the closed ternary ideal $vG$ of $E$, where $G$ is ternary isomorphic to $vG$ via $v$. It follows that $\cA_G\cong\cB_{vG}$ via $\psi\colon\AB{g,g'}\mapsto\AB{vg,vg'}$. Composing $x\in\sB^a(\cA_G,G)$ to $vx\psi^{-1}\in\sB^a(\cB_{vG},vG)$ establishes a ternary (and strict!) isomorphism $M(G)\rightarrow M(vG)$ which, of course, passes on to an isomorphism of the quotients $Q(G)\cong Q(vG)$.

So, summarizing, the Busby invariant of an extension $E$ of $F$ by $G$ is $0$ if and only if $E$ is isomorphic to the external direct sum $F\oplus G$ and the ternary homomorphism $u$ and $v$ correspond, under this isomorphism, to the canonical maps making $F\oplus G$ an extension. So, the trivial Busby invariant corresponds exactly to the super-trivial extension. As a corollary we recover the statement, well-known for extensions of \nbd{C^*}algebras, that extensions by a Hilbert module $G$ with trivial corona space $Q(G)=\zero$ (that is, $M(G)=G$) are the super-trivial ones.
\end{enumerate}
\eobn
\eemp

\vspace{-2ex}
\bemp[``Adding'' extensions.~] \label{Badd}
Homomorphisms can be added in various ways involving direct sums. For adding the Busby invariants $\tau_i\colon\cC\rightarrow Q(\cA)$ ($i=1,2$)  of two extensions of $\cC$ by $\cA$ (both fixed), one starts with the direct sum $\tau_1\oplus\tau_2\colon\cC\rightarrow Q(\cA)\oplus Q(\cA)=Q(\cA\oplus\cA)$ that sends $c$ to $\tau_1(c)\oplus\tau_2(c)$. Without any further assumptions on $\cA$, this is the Busby invariant of an extension of $\cC$ by $\cA\oplus\cA$. One might think about the case when $\cA\oplus\cA$ is isomorphic to $\cA$. But this would, for instance, exclude simple \nbd{C^*}algebras. But, $\cA\oplus\cA$ is the diagonal subalgebra of $M_2(\cA)$, and one may view $\tau_1\oplus\tau_2$ as homomorphism from $\cC$ into $M_2(\cA)$. (At least, $M_2(\cA)$ is simple, if $\cA$ is; more precisely, the centers of $M_2(\cA)$ and $\cA$ coincide.) It is not so rare (for instance, if $\cA$ is \it{stable}) that $\cA\cong M_2(\cA)$. So, at least in that case, we may add Busby invariants up to automorphisms of $M_2(\cA)\cong\cA$, and passing to the corresponding (equivalence class of) extension(s), we have a ``sum'' operation among extensions of $\cC$ by $\cA$ (actually, a monoid). The restriction $M_2(\cA)\cong\cA$ is not as harsh as it might sound. Anyway, we know that the extensions are super-trivial, if $\cA$ is unital. And if $\cA$ is nonunital, even if $\cA$ should not be isomorphic to $M_2(\cA)$, we may pass to the \hl{stabilization} $\sK\otimes\cA$ (where $\sK:=\sK(\K)$ for an infinite-dimensional separable Hilbert space $\K$), which preserves many of the properties of $\cA$ (properties up to \it{Morita equivalence}, to be precise).

We can do the same for extensions of $F$ by $G$; and, again, we do it via the corresponding blockwise extensions of the linking algebras by Theorem \ref{extthm}. In order to proceed that way, we have to hypothesize that $\sK\rtMatrix{\cA_G\\G}\cong M_2\family{\sK\rtMatrix{\cA_G\\G}}$. And in order to speak about blockwise extensions, we have to interpret the latter as reduced linking algebra via $M_2\family{\sK\rtMatrix{\cA_G\\G}}=\sK\rtMatrix{M_2(\cA_G)\\M_2(G)}$. Here, $M_2(G)$ is a Hilbert \nbd{M_2(\cA)}module in an obvious fashion, full over $M_2(\cA_G)$, so that $\sK\rtMatrix{M_2(\cA_G)\\M_2(G)}$ is, indeed, its reduced linking algebra. (It is important to note that the blockwise structure of the matrix $M_2\family{\sK\rtMatrix{\cA_G\\G}}$ has \bf{absolutely nothing} to do with the blockwise structure of the reduced linking algebra $\sK\rtMatrix{M_2(\cA_G)\\M_2(G)}$! The former is a \nbd{2\times 2}matrix of reduced linking algebras; the latter is a reduced linking algebra whose blocks consist of \nbd{2\times 2}matrices.)

If we have $\sK\rtMatrix{\cA_G\\G}\cong M_2\family{\sK\rtMatrix{\cA_G\\G}}=\sK\rtMatrix{M_2(\cA_G)\\M_2(G)}$, then this isomorphism turns over to the multiplier algebras $\sB^a\rtMatrix{\cA_G\\G}=M\family{\sK\rtMatrix{\cA_G\\G}}\cong M\family{M_2\family{\sK\rtMatrix{\cA_G\\G}}}=M\family{\sK\rtMatrix{M_2(\cA_G)\\M_2(G)}}=\sB^a\rtMatrix{M_2(\cA_G)\\M_2(G)}$ and the corona algebras $Q\family{\sK\rtMatrix{\cA_G\\G}}\cong Q\family{\sK\rtMatrix{M_2(Q(\cA_G))\\M_2(Q(G))}}$. The direct sum $Q\family{\sK\rtMatrix{\cA_G\\G}}\oplus Q\family{\sK\rtMatrix{\cA_G\\G}}=:\rtMatrix{Q_{1,1}&Q_{1,2}\\Q_{2,1}&Q_{2,2}}\oplus\rtMatrix{Q_{1,1}&Q_{1,2}\\Q_{2,1}&Q_{2,2}}$ sits in $Q\family{\sK\rtMatrix{M_2(Q(\cA_G))\\M_2(Q(G))}}$ as
\beqn{
\Matrix{\rtMatrix{Q_{1,1}&~\\~&~}&\rtMatrix{Q_{1,2}&~\\~&~}\\\rtMatrix{Q_{2,1}&~\\~&~}&\rtMatrix{Q_{2,2}&~\\~&~}}
~\oplus~
\Matrix{\rtMatrix{~&~\\~&Q_{1,1}}&\rtMatrix{~&~\\~&Q_{1,2}}\\\rtMatrix{~&~\\~&Q_{2,1}}&\rtMatrix{~&~\\~&Q_{2,2}}}.
}\eeqn
(Likewise for $\sK\rtMatrix{\cA_G\\G}$ and $\sB^a\rtMatrix{\cA_G\\G}$, but for the Busby invariant we need this only for $Q$.) So it is clear, how we have to add two blockwise Busby invariants $\sT_i=\rtMatrix{\tau_i&T_i^*\\T_i&\vt_i}$ to obtain their ``sum''
\beqn{
\sT
~=~
\sT_1\oplus\sT_2
\colon
\sK\rsMatrix{\cC_F\\F}
~\longrightarrow~
Q\family{\sK\rsMatrix{\cA_G\\G}}\oplus Q\family{\sK\rsMatrix{\cA_G\\G}}
~\subset~
Q\family{\sK\rsMatrix{M_2(Q(\cA_G))\\M_2(Q(G))}}
~\cong~
Q\family{\sK\rsMatrix{\cA_G\\G}}.
}\eeqn
The \nbd{21}corner of $\sT$ is a map $T\colon F\rightarrow Q(G)\oplus Q(G)\subset M_2(Q(G))\cong Q(G)$, where the last isomorphism is the ternary isomorphism obtained as the \nbd{21}corner of the blockwise isomorphism $Q\family{\sK\rtMatrix{\cA_G\\G}}\cong Q\family{\sK\rtMatrix{M_2(Q(\cA_G))\\M_2(Q(G))}}$. This map $T$ is the \hl{sum} (up to suitable equivalence) of the Busby invariants $T_1$ and $T_2$.

The sum of the Busby invariants may be defined directly by the preceding line, as soon as we have a ternary isomorphism $Q(G)\cong M_2(Q(G))$. However, only the derivation via the corresponding block-wise extensions makes such a definition not appear to be \it{ad hoc} and, more importantly, makes available in a blink of an eye all the related structures (like Ext-semigroups) and the results on them, well known for extensions of \nbd{C^*}algebras. We briefly address the question, when such an isomorphism $Q(G)\cong M_2(Q(G))$ exists.

Necessarily, $\cA_G\cong M_2(\cA_G)$ (and also $\sK(G)\cong M_2(\sK(G))$), because by Theorem \ref{homothm} a ternary isomorphism lifts to a blockwise isomorphism between the reduced linking algebras. However, also an opposite statement is true when we start with $\cA_G\cong M_2(\cA_G)$ (or $\sK(G)\cong M_2(\sK(G))$). Indeed, if we fix an isomorphism $\vp\colon\cA_G\rightarrow M_2(\cA_G)$, then $ga\mapsto g\odot a\mapsto g\odot\vp(a)\mapsto\rtMatrix{g\odot\vp(a)_{1,1}&g\odot\vp(a)_{1,2}\\g\odot\vp(a)_{2,1}&g\odot\vp(a)_{2,2}}\mapsto\rtMatrix{g\vp(a)_{1,1}&g\vp(a)_{1,2}\\g\vp(a)_{2,1}&g\vp(a)_{2,2}}$ is a ternary isomorphism $G=G\odot\cA\cong G\odot M_2(\cA)\cong M_2(G\odot\cA)=M_2(G)$. (The first and last equality refer to the canonical (unitary) identification of a Hilbert \nbd{\cB}module $E$ with $E\odot\cB$, and also the isomorphism $E\odot M_2(\cB)\cong M_2(E\odot\cB)$ is unitary. (Exercise: Find the correct action of $\cB$ on $M_2(\cB)$ that makes tensor product $E\odot M_2(\cB)$ over $\cB$ isomorphic to $M_2(E\odot\cB)$.) Only the step from $G\odot\cA$ to $G\odot M_2(\cA)$ is a ternary isomorphism that goes beyond unitaries; see, for instance, \cite{Ske06c} for more details.) Note that we may replace $\cA$ with $\cA_G$. (The statement starting from an isomorphism $\sK(G)\cong M_2(\sK(G))$ follows by symmetry of the reduced linking algebra under the exchange $G\leftrightarrow G^*$.) We just proved the following result:

\bpropn
~~~~~~~~~$\cA_G\cong M_2(\cA_G)$ ~~$\Longleftrightarrow$~~ $G\cong M_2(G)$ ~~$\Longleftrightarrow$~~ $\sK(G)\cong M_2(\sK(G))$.
\epropn

\noindent
So addition of extensions of $F_\cC$ by $G_\cA$ makes sense under exactly the same conditions (on $\cA_G$, only) under which addition of extensions of $\cC_F$ by $\cA_G$ makes sense.
\eemp

\brem \label{extrem}
Here, we stop with our brief on extensions of Hilbert modules. We do not make any attempt to compare in detail the sketchy insights of this section with the work of Bakic and Guljas \cite{BaGu02,BaGu04} and the forthcoming papers (also by others) in more detail than the few references that popped already up throughout this section.

Surely their work is much more complete. But we think that adding Busby invariants might be new in these notes. Of course, this opens up to look at all the known structures (like Ext (semi)groups), well-known for \nbd{C^*}algebras, also for Hilbert modules. (Needless to say that we recommend to do this in terms of the reduced linking algebras, and, then, reduce it to an intrinsic formulation in terms of ternary morphisms.) We also do not know if the insights in Observation \ref{Bob3} under \ref{Badd} did already exist or are new.

Our scope was just to illustrate that many results and -- almost more importantly -- also definitions about extensions of Hilbert modules drop out by looking at well-known things for reduced linking algebras. (One problem that occurred also in \cite{BaGu03} and that we could not really get rid of, is the fact that the multiplier space of a full Hilbert module is, usually, not full over the multiplier algebra; this is an inconvenience which remains.)

It might be interesting to have a look how this section fits into the context of \it{triple extensions} as they occur in the context of \it{noncommutative Shilov boundary}. (See Blecher and Le Merdy \cite[Section 8.3]{BlLMe04} and the references mentioned in \cite[Section 8.7]{BlLMe04}. They seem to fit in the context of \it{essential extensions} in \cite{BaGu04} and the connections with the work of Hamana (see \cite{BlLMe04}) could also have been discussed in \cite{BaGu04}.) We omit also that.
\erem

\newpage

\section{Questions} \label{qSEC}

We close with three (sets of) questions, that follow our philosophy of motivating notions for Hilbert modules and proving results about them, by looking at known notions and their results for the associated reduced linking algebras.

\lf
\bemp \label{Q1}
We mentioned that a split extension $\cB$ is, roughly, a conditional expectation $\E$ ($=s\circ\vp$) onto a subalgebra ($s(\cC)$) such that $\ker\E$ is an ideal. (See below, for conditional expectation.) Forgetting about the condition that $\ker\E$ is an ideal (that arises from the special situation in a split extension, where $\E$ has to be a homomorphism), it is natural to ask what is a conditional expectation from a Hilbert module onto a closed ternary subspace. (Recall from Lemma \ref{TsSlem} that this is a weaker condition than requiring a closed submodule; submodules are not what corresponds to subalgebras but to right ideals; ternary subspaces are what corresponds to subalgebras.) So, accepting that a \it{conditional expectation} of a Hilbert module onto a closed ternary subspace should be equivalently defined as (co)restriction of a blockwise conditional expectation between the reduced linking algebras, what are its intrinsic properties (as maps between modules) that guarantee to make this happen?

Recall that a \hl{conditional expectation} $\E$ of \nbd{C^*}algebras is a linear idempotent from a \nbd{C^*}al\-gebra $\cA$ \bf{onto} a \nbd{C^*}subalgebra $\cB$ fulfilling one the following two equivalent (see Tomiyama \cite{Tom57} and Takesaki \cite[Theorem III.3.4]{Tak02}) conditions:
\begin{enumerate}
\item[(i)]
$\E$ is a contraction. (More precisely, $\norm{\E}=1$ unless $\cB=\zero$.)

\item[(ii)]
$\E$ is a positive \nbd{\cB}bimodule map.
\end{enumerate}
So, assuming we have a blockwise conditional expectation $\E\colon\sK\rtMatrix{\cB_E\\E}\rightarrow\sK\rtMatrix{\cC_F\\F}$, the (co)re\-strictions $\E_{i,j}$ to each \nbd{ij}corner are clearly idempotent contractions. So, the (co)re\-stric\-tions $\E_{i,i}$ to the diagonal entries are conditional expectations in their own right. A whole bunch of algebraic conditions follow, when we write down what the condition to be a bimodule map means for the restriction to products from different corners. The two questions we ask here are: For an idempotent from $E$ onto a closed ternary subspace $F$, is it enough to check only the ``ternary'' condition $\E(\E(x)\AB{y,\E(z)})=\E(x)\AB{\E(y),\E(z)}$ to guarantee an extension to a blockwise conditional expectation? (For instance, we are not sure if this allows to conclude that this map is a contraction, before actually having the blockwise extension.) Is it possibly even enough, to require that the idempotent is a contraction?

(We should keep in mind that it might be a good idea to pass to von Neumann algebras and von Neumann modules. Also in proving equivalence of (i) and (ii), passing to the biduals is an essential step in either direction. Let us also mention that Hahn-Banach type extension of off-diagonal (complete) contractions to (completely) positive contractions are not (directly) applicable; such extensions are limited to \it{injective} codomains. We also should keep in mind that a conditional expectation is a CP-map; as such it has a GNS-construction -- a GNS-construction with special properties. Being a blockwise CP-map between linking algebras, relates the present question to the following.)
\eemp

\lf
\bemp \label{Q2}
\hl{Semisplit extensions} are like split extensions, just that the \it{splitting} $s$ is allowed to be a CP-map instead of a homomorphism. This leads directly to the question, what is a CP-map between Hilbert modules. In Skeide and Sumesh \cite{SkSu14} we proposed -- surprise! -- that a map between (full) Hilbert modules is \hl{CP} if it extends to a blockwise CP-map between the linking algebras, to be precise between the \it{extended} linking algebras, that is \it{strict} on its \nbd{22}corner. (The \hl{extended linking algebra} of a Hilbert module $E$ is $\rtMatrix{\cB&E^*\\E&\sB^a(E)}\subset\sB^a\rtMatrix{\cB\\E}$. This enlargement is unavoidable for the codomain, because, as Stinespring type constructions show, the \nbd{22}corner involves amplifications $\sK(E)\odot\id$, which, in general, are no longer compact operators. Doing this consequently (and required once we wish to compose CP-maps), also in the domain we replace $\sK(E)$ with $\sB^a(E)$; and being \it{strict} is just the natural (and indispensable) compatibility condition (which, under certain nondegeneracy conditions, also is fulfilled automatically).) There are authors who propose Asadi's \nbd{\tau}maps as CP-maps between Hilbert modules. (For a linear map $\tau\colon\cB\rightarrow\cC$, a map $T\colon E_\cB\rightarrow F_\cC$ between Hilbert modules is a \nbd{\tau}map, if $\AB{Tx,Tx'}=\tau(\AB{x,x'})$. They have been proposed for CP-maps $\tau$ and $\cC=\sB(H), F=\sB(H,H')$ in Asadi \cite{Asa09}. Bhat, Ramesh, and Sumesh \cite{BRS12} (correctly) proved a Stinespring type theorem suggested by Asadi, also removing a \it{trivializing} condition from the hypotheses in \cite{Asa09}. Skeide \cite{Ske12a} presented a half-a-page proof using module language, for general $\cC$ and $F$. A discussion of \nbd{\tau}maps which is quite exhaustive in many senses followed in \cite{SkSu14}. They also show that bounded $\tau$ is CP on $\cB_E$, automatically.)

We do not think that the restriction for CP-maps between modules to be also \nbd{\tau}maps is justified. Following our philosophy that a semisplit extension of (full, for safety) Hilbert modules corresponds/has to do with blockwise semisplit extensions of the (extended) linking algebras, our definition of CP-map promises to be the more suitable. (Why should we restrict the blockwise CP-maps making a blockwise extension semisplit to a subset of the blockwise CP-maps?) Independently, one may examine also the subclass of semisplit extensions where the splitting is required to be a \nbd{\tau}map, and analyze which special properties they share.
\eemp

\lf
\bemp \label{Q3}
Ideals in a \nbd{C^*}algebra $\cB$ are linear subspaces $I$ such that $\cB I\cB\subset I$. \hl{Hereditary subalgebras} can be equivalently characterized as \nbd{C^*}subalgebras $\cC$ such that $\cC\cB\cC\subset\cC$. A pragmatic (that is, blockwise) definition for a closed ternary subspace $F$ of a Hilbert \nbd{\cB}module $E$ to be \hl{linking hereditary} would be the requirement that the reduced linking algebra of $F$ is hereditary in the reduced linking algebra of $E$. A candidate for an intrinsic definition would be to say $F$ is \hl{ternary hereditary} in $E$ if $F\AB{E,F}\subset F$. Clearly, the former implies the latter. But, the former also implies that $\cC_F$ is hereditary in $\cB_E$ and that $\sK(F)$ is hereditary in $\sK(E)$. However, we presently do not see a reason, why a hereditary ternary subspace would imply that $\cB_F$ and $\sK(F)$ are hereditary, too. On the other hand, if $\cB_F$ or $\sK(F)$ is hereditary, then
\bmun{
F\AB{E,F}
~\subset~
\cls\sK(F)F\AB{E,F}\cB_F
~=~
\cls F\AB{F,F}\AB{E,F}\AB{F,F}
\\
~=~
\cls(FF^*)(FE^*)(FF^*)F
~\subset~
F,
}\emun
because $\AB{F,F}\AB{E,F}\AB{F,F}\subset\cB_F$ or because $(FF^*)(FE^*)(FF^*)\subset\sK(F)$. Again we do not know if one of $\cB_F$ and $\sK(F)$ being hereditary would imply the other one being hereditary, too. But at least we get: A closed ternary subspace $F$ is linking hereditary if and only of both $\cB_F$ and $\sK(F)$ are hereditary.

Clearly, closed submodules are linking hereditary. (Then $\cB_F$ is even an ideal in $\cB_E$ and, clearly $\sK(F)$ is hereditary in $\sK(E)$.) So, each of the two conditions in \ref{TI=TI} separately imply that $F$ is linking hereditary. (Indeed, the second condition $F\AB{E,E}\subset F$ implies that $F$ is Hilbert submodule of $E$.) This is like the obvious statement that closed right or left ideals in \nbd{C^*}algebras are hereditary subalgebras.

We leave open the question whether or not ternary hereditary implies linking hereditary. (Of course, if not, then we think the ``right'' definition of hereditary subspaces is linking hereditary.) What we call a closed ternary hereditary subspace, occurred as \it{inner ideal} of a TRO in \cite{BlNe07}.

There is another much more interesting (and probably more far reaching) question, namely: The usual definition of hereditary subalgebra is that the positive elements of $\cC$ form a \hl{hereditary subcone} of the positive elements of $\cB$, that is, $0\le b\le c$ for $b\in\cB$ and $c\in\cC$ implies $b\in\cC$. (We preferred to start with the equivalent condition $\cC\cB\cC\subset\cC$, because this condition immediately suggests analogue properties for Hilbert modules.) Is there a similar possibility to describe exactly the same (linking) hereditary subspaces of a Hilbert module in terms of positive cones?

We may introduce a pre-order on a Hilbert module $E$ by saying $x\le y$ if $\AB{x,x}\le\AB{y,y}$. (Kolarec \cite{Kol10} proposed $\sqrt{\AB{x,x}}\le\sqrt{\AB{y,y}}$. Clearly, the former implies the latter but not \it{vice versa}.) We don't know, if linking hereditary subspaces $F$ are \hl{right hereditary} in the sense that $x\le y$ for all $y\in F$ implies $x\in F$. (We think rather not.) But, if $\cC$ is a hereditary subalgebra of $\cB$, then $F:=\cls E\cC$ is a Hilbert \nbd{\cC}module and right hereditary in $E$. Something similar goes for the \hl{left hereditary} subspace $\cls\cA E$ for a hereditary subalgebra $\cA$ of $\sK(E)$, but then rather for the Hilbert \nbd{\sK(E)}modules $E^*$.

There are papers about positivity in Hilbert modules. (See, for instance, Blecher and Neal \cite{BlNe07} and the references therein.) But, we should note that positivity in a TRO $E$ is referring to positivity in $E\cap E^*$ in the \nbd{C^*}algebra in which the TRO sits. If this \nbd{C^*}algebra is the linking algebra, then this intersection is $\zero$. So, not only is this positivity not compatible with our point of view to take intuition from the linking algebras, but is also depends on the choice of an embedding of $E$ as a TRO into a \nbd{C^*}algebra. What we ask for is a structure that is intrinsic to the (full) Hilbert module $E$, a structure that is present for every Hilbert module (and tells enough about it to determine most of it, to the same extent as the positive cone of a \nbd{C^*}algebra tells about the structure of the \nbd{C^*}algebra), and a structure that allows to recover what we defined to be a linking hereditary subspace.
\eemp



\newpage

 \noindent
 \bf{Acknowledgment:} I wish to express a big thank you to the referee of this version.

\setlength{\baselineskip}{2.5ex}

\newcommand{\Swap}[2]{#2#1}\newcommand{\Sort}[1]{}
\providecommand{\bysame}{\leavevmode\hbox to3em{\hrulefill}\thinspace}
\providecommand{\MR}{\relax\ifhmode\unskip\space\fi MR }
\providecommand{\MRhref}[2]{%
  \href{http://www.ams.org/mathscinet-getitem?mr=#1}{#2}
}
\providecommand{\href}[2]{#2}

\lf\noindent
Michael Skeide:
{\small\itshape Dipartimento di Economia, Università degli Studi del Molise, Via de Sanctis, 86100 Campobasso, Italy, E-mail: \href{mailto:skeide@unimol.it}{\tt{skeide@unimol.it}}}\\
{\small{\itshape Homepage: \url{http://web.unimol.it/skeide/}}}


\end{document}